\documentclass{amsart}
\usepackage{amsmath, amssymb}
\usepackage[pdftex]{graphicx, xcolor}  
\usepackage{latexsym}

\theoremstyle{definition}

\newtheorem{remark}{Remark}

\newcommand{\red}[1]{\textcolor{red}{#1}}

\renewcommand{\theenumi}{\roman{enumi}}

\title[Numerical analysis of evolving spirals] 
      {Numerical analysis comparing ODE approach and level set method
for evolving spirals by crystalline eikonal-curvature flow}

\author[Tetsuya Ishiwata and Takeshi Ohtsuka]{}

\subjclass[2010]{Primary: 34A34, 53C44; Secondary: 53A04.}
 \keywords{Crystalline eikonal-curvature flow, 
 Evolution of a polygonal spiral,
 Level set method.}

 \email{tisiwata@shibaura-it.ac.jp}
 \email{tohtsuka@gunma-u.ac.jp}

\thanks{The first author is
partly supported by JSPS KAKENHI Grant Number
15H03632 and 16H03953.}

\begin{document}
\maketitle

\centerline{\scshape Tetsuya Ishiwata}
\medskip
{\footnotesize
 \centerline{Department of Mathematical Sciences, 
Shibaura Institute of Technology}
   \centerline{Fukasakuk 309, Minuma-ku, Saitama 337-8570, Japan}
} 

\medskip

\centerline{\scshape Takeshi Ohtsuka} 
\medskip
{\footnotesize
 \centerline{Division of Pure and Applied Science,
 Faculty of Science and Technology, Gunma University}
   \centerline{Aramaki-machi 4-2, Maebashi, 371-8510 Gunma, Japan}
}

\bigskip


\begin{abstract}
 In this paper, the evolution of a polygonal spiral curve
 by the crystalline curvature flow with a pinned center 
 is considered with two view points, discrete model 
 consist of an ODE system of facet lengths 
 and a level set method.
 We investigate the difference of these models numerically
 by calculating the area of the region enclosed by these spiral curves.
 The area difference is calculated by the normalized
 $L^1$ norm of the difference of step-like functions
 which are branches of $\arg x$ whose discontinuities
 are only on the spirals.
 We find the differences of the numerical results
 considered in this paper are very small 
 even though the evolution laws of these models
 around the center and the farthest facet
 are slightly different.
\end{abstract}

\section{Introduction}


The crystalline curvature of a curve $\Gamma$, 
which is denoted by $H_\gamma$,
is defined by the changing ratio of an anisotropic surface
energy functional
\[
 E_\gamma (\Gamma) = \int_\Gamma \gamma ({\bf n}) d \sigma
\]
for a singular
 density function 
$\gamma \colon \mathbb{R}^2 \to [0,\infty)$
with respect to the volume of a region enclosed by $\Gamma$,
where ${\bf n}$ is a continuous unit normal vector field of
$\Gamma$ and $d \sigma$ is the line element.
Here, singular means that the Wulff shape
\[
 \mathcal{W}_\gamma = \{ p \in \mathbb{R}^2 ; \ 
 p \cdot q \le \gamma (q) \ \mbox{for} \ q \in \mathbb{S}^1 \},
\]
which satisfies $H_\gamma = 1$ on $\partial \mathcal{W}_\gamma$,
is a convex polygon.
See \cite{Gurtin:1993} for details of the crystalline curvature.
Such a singular energy expresses the surface energy
of the polygonal structure of interfaces like as 
crystal surface.
The typical example of $\gamma$ is $\ell^1$ norm.
For describing general settings,
we here assume that
\begin{enumerate}
 \renewcommand{\theenumi}{A\arabic{enumi}}
 \item \label{assum: convexity}
       $\gamma$ is convex,
 \item \label{assum: phd1}
       $\gamma$ is positively homogeneous of degree 1, i.e.,
       $\gamma (\lambda p) = \lambda \gamma (p)$
       for $p \in \mathbb{R}^2$ and $\lambda > 0$,
 \item \label{assum: positivity}
       $\gamma > 0$ on $\mathbb{S}^1$
 \item \label{assum: crystalline}
       $\gamma$ is piecewise linear.
\end{enumerate}
Note that \eqref{assum: phd1} is 
for the level set formulation of curves
mentioned later.
Moreover,
(\ref{assum: crystalline})
is a sufficient condition 
to the singularity of $\mathcal{W}_\gamma$
for the crystalline curvature,
since
$\mathcal{W}_\gamma = \{ p \in \mathbb{R}^2 ; \ \gamma^\circ (p) \le 1 \}$
and $(\gamma^\circ)^\circ = \gamma$
if $\gamma$ is convex, where
$\gamma^\circ (p) := \sup \{ p \cdot q ; \ \gamma (q) \le 1 \}$
is a support function of $\gamma$.
See \cite{Rockafellar:book} for details of 
the properties of $\gamma$ and $\gamma^\circ$.

In this paper we consider the evolution of 
a convex polygonal spiral by 
\begin{equation}
 \label{geo mcf}
  \beta V_\gamma = U - \rho_c H_\gamma
  \quad \mbox{on} \ \Gamma_t,
\end{equation}
where $V_\gamma$ is an anisotropic normal velocity
under the Finsler metric defined by 
$\mathrm{dist}_\gamma (x,y) = \gamma^\circ (x-y)$,
and $U > 0$ and $\rho_c > 0$ are assumed to be constants.
(Note that we \emph{do not assume} the symmetricity
of this metric.)
For this evolution of a pinned spiral, the authors of this paper 
introduce a discrete model by an ODE system of the facet lengths
in \cite{IO:DCDS-B},
due to the idea of \cite{AG:1989ARMA, Taylor1991, Ishiwata2014},
see also \cite{Gurtin:1993} for details.

On the other hand, Tsai, Giga and the second author 
\cite{Ohtsuka:2003wi, OTG:2015JSC} introduced 
a level set formulation for evolving spirals with fixed centers.
According to their formulation, an evolving spiral curve
with a fixed center at the origin
is given as 
\[
 \Gamma_L (t) = \{ x ; \ u(t,x) - \theta (x) \equiv 0 \mod
 2 \pi \mathbb{Z} \}, \quad
 {\bf n} = - \frac{\nabla (u - \theta)}{|\nabla (u - \theta)|}
\]
with an auxiliary function $u(t,x)$ and a pre-defined
multivalued function $\theta (x) = \arg x$.
Then, $V_\gamma$ and $H_\gamma$ are interpreted as
\[
 V_\gamma = \frac{u_t}{\gamma (- \nabla (u - \theta))}, \quad
 H_\gamma = - \mathrm{div} \{ \xi (- \nabla (u - \theta)) \},
\]
where $\xi = D \gamma$.
Hence, we obtain the level set equation
for \eqref{geo mcf}
of the form
\[
  \tilde{\beta} (\nabla (u - \theta)) u_t 
  - \tilde{\gamma} (\nabla (u - \theta))
  \left[
   \mathrm{div} \{ \tilde{\xi} (\nabla (u - \theta)) \}
   + U
  \right] = 0,
\]
where 
$\tilde{\beta} (p) = \beta (-p)$,
$\tilde{\gamma} (p) = \gamma (-p)$,
and $\tilde{\xi} (p) = \xi (-p)$.

The aim of this paper is to show the numerical difference
between the spirals calculated by the discrete model
due to \cite{IO:DCDS-B} and the level set method
due to \cite{OTG:2015JSC}.
To measure the difference between these spirals,
we calculate the area of the region enclosed
by their spirals.
It is established by calculating
\[
 \mathcal{D} (t)
 = \frac{1}{|W|} \int_W
 \frac{\theta_D (t,x) - \theta_L (t,x)}{2 \pi} dx,
\]
where $\theta_D$ and $\theta_L$ are branches of
$\theta$ whose discontinuities are only on
the spiral curves 
$\Gamma_D (t) = \sum_{j=0}^k L_j (t)$
obtained by the discrete algorithm and
$\Gamma_L (t)$
by the level set method, respectively.
A practical way to construct $\theta_L$
from solution $u$ of the level set equation
is provided in \cite{OTG:2015JSC}.
Thus, we shall give a way to construct $\theta_D$
in \S \ref{sec: diff function}.
Note that the discrete model in \cite{IO:DCDS-B} is
constructed from $\mathcal{W}_\gamma$,
we shall give a way to construct $\gamma$ from $\gamma^\circ$
in \S \ref{sec: how to define density} to obtain
the level set equation corresponding to the discrete model.

\section{Models}

In this section, we recall the discrete model
due to \cite{IO:DCDS-B} and 
the level set method due to \cite{OTG:2015JSC}.
To compare the evolving spiral curves from these models,
we 
have to give a Wulff shape $\mathcal{W}_\gamma$
for the discrete model 
and corresponding surface energy density $\gamma$
for the level set method.
In this section, we consider the situation 
$\mathcal{W}_\gamma$ and corresponding $\gamma$
are already given.
A practical way to obtain $\gamma$ from $\mathcal{W}_\gamma$
will be discussed in \S \ref{sec: how to define density}.
We briefly review mathematical results on these models.

\subsection{Discrete model}
\label{sec: ODE model}

We recall the ODE model by \cite{IO:DCDS-B}.

We first prepare some notations for $\mathcal{W}_\gamma$.
Let $\mathcal{W}_\gamma$ be a $N_\gamma$ sided 
convex polygon.
The $j$-th facet of $\mathcal{W}_\gamma$
has an outer unit normal vector $\mathcal{N}_j$
with angle $\varphi_j$ for $j=0, 1, 2, \ldots, N_\gamma - 1$.
Set the unit tangential vector ${\bf T}_j$
of the $j$-th facet as well as the definition
of the Frenet frame, i.e.,
\[
 {\bf N}_j = (\cos \varphi_j, \sin \varphi_j), \quad
 {\bf T}_j = (\sin \varphi_j, - \cos \varphi_j).
\]
We assume the followings for
expressing the convexity of $\mathcal{W}_\gamma$.
\begin{enumerate}
 \renewcommand{\theenumi}{W\arabic{enumi}}
 \item $\varphi_0 < \varphi_1 < \varphi_2 < \cdots < \varphi_{N_\gamma - 1}
       < \varphi_0 + 2 \pi$.
 \item $\varphi_j < \varphi_{j+1} < \varphi_j + \pi$
       for $j=0,1,2,\ldots, N_\gamma - 1$.
\end{enumerate}
Note that $\varphi_{N_\gamma} = \varphi_0$.
We denote the length of the $j$-th facet of $\mathcal{W}_\gamma$
by $\ell_j > 0$.

\begin{figure}[htbp]
 \begin{center}
  \includegraphics[scale=0.9]{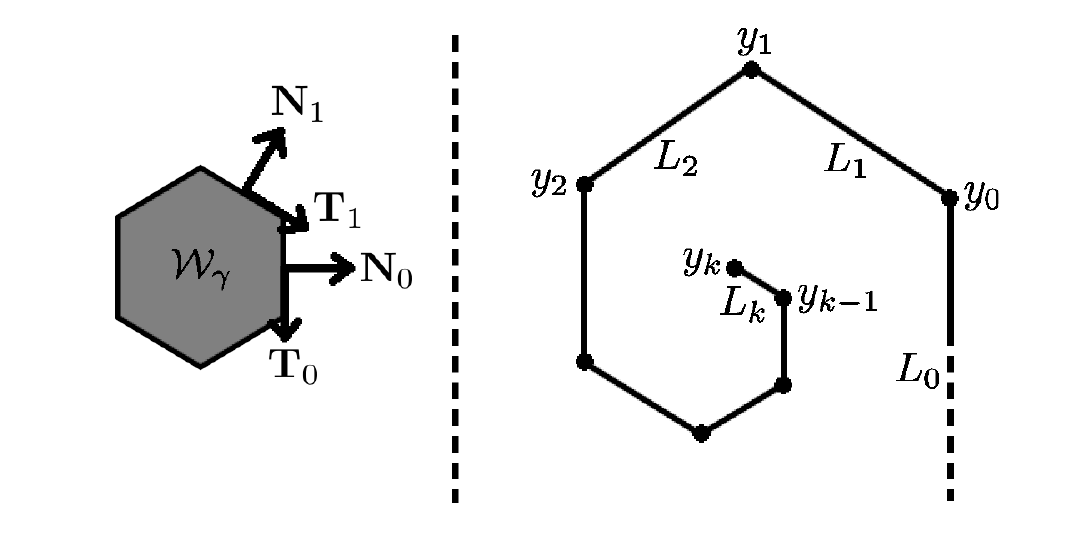}
  \caption{Description of $\Gamma_D = \bigcup_{j=0}^k L_j (t)$.
  Note that the variable $t$ of $L_j$ and $y_j$ is omitted 
  in the above figure for the simplicity.
  }
  \label{facet_vertices}
 \end{center}
\end{figure}

We next prepare the notation of an evolving 
polygonal spiral. 
We denote an evolving polygonal spiral curve by \eqref{geo mcf}
by $\Gamma_{D} (t) = \bigcup_{j=0}^k L_j (t)$.
According to \cite{IO:DCDS-B},
we here consider the evolution of a positive
convex polygonal spiral.
Assume that the $j$-th facet $L_j (t)$ is given as
\[
 L_j (t) = 
 \left\{
 \begin{array}{ll}
  \{ \lambda y_j (t) + (1 - \lambda) y_{j-1} (t) ; \
   \lambda \in [0,1] \} & \mbox{for} \ j=k, k-1, \ldots, 1 \\
  \{ y_0 (t) + \lambda {\bf T}_0 ; \ \lambda > 0 \} &
   \mbox{if} \ j = 0
 \end{array}
 \right.
\]
with vertices $y_j (t)$ ($j = 0,1,2, \ldots, k-1$)
and the center $y_k (t) = O$.
Assume that 
\[
 {\bf T}_j = \frac{y_{j-1} (t) - y_j (t)}{|y_{j-1} (t) - y_j (t)|}.
\]
We have extended the number $j$ of ${\bf T}_j$ 
from $j=0.1,2, \ldots, N_\gamma - 1$ to $\mathbb{Z}$;
let ${\bf T}_{j + nN_\gamma} = {\bf T}_j$
for $j=0,1, 2, \ldots, N_\gamma -1$
and $n \in \mathbb{Z}$.
Then, the evolution of $\Gamma_D (t)$
by \eqref{geo mcf} with fixed center $y_k (t) = O$
is expressed by an ODE system for 
$d_j (t) = |y_j (t) - y_{j-1} (t)|$ of the form
\begin{align}
 \label{ODE system1}
 & \dot{d}_k = c^-_k \left(U - \frac{\rho_c \ell_{k-1}}{d_{k-1}} \right), \\
 \label{ODE system2}
 & \left\{
 \begin{aligned}
  \dot{d}_{k-1} 
  & = - b_{k-1} \left( U - \frac{\rho_c \ell_{k-1}}{d_{k-1}} \right)
  + c^-_{k-1}
  \left( U - \frac{\rho_c \ell_{k-2}}{d_{k-2}} \right), \\
  \dot{d}_j 
  & = - b_j \left( U - \frac{\rho_c \ell_j}{d_j} \right)
  + c^+_j
  \left( U - \frac{\rho_c \ell_{j+1}}{d_{j+1}} \right)
  + c^-_j
  \left( U - \frac{\rho_c \ell_{j-1}}{d_{j-1}} \right) \\
  & \hspace{6cm} \mbox{for} \ j=2, 3, \ldots, k-2, \\
  \dot{d}_1 
  & = - b_1 \left( U - \frac{\rho_c \ell_1}{d_1} \right)
  + c^+_1
  \left( U - \frac{\rho_c \ell_2}{d_2} \right)
  + c^-_1 U,
 \end{aligned}
 \right.
\end{align}
where $b_j \in \mathbb{R}$ and $c^\pm_j > 0$
are numerical constants defined by
\[
 b_j = \frac{1}{\beta_j}
 \left( \frac{1}{\tan(\varphi_{j+1} - \varphi_j)} 
 + \frac{1}{\tan(\varphi_j - \varphi_{j-1})} \right), \quad
 c_j^\pm = \pm \frac{1}{\beta_{j \pm 1} \sin (\varphi_{j \pm 1} - \varphi_j)}
\]
and $\beta_j = \beta ({\bf N}_j)$.
Tracking the evolution of $\Gamma_D (t)$ is established 
by drawing $\Gamma_D (t)$ with setting
\[
 y_k (t) = O, \quad
 y_{j-1} (t) = y_j (t) + d_j (t) {\bf T}_j \
 \mbox{for} \ j=k, k-1, k-2, \ldots, 1.
\]
See Figure \ref{facet_vertices} for 
details of $\Gamma_D (t)$ described with
the above notations.

In this paper, we give an initial curve as
$k=1$ with $d_1 (0) = 0$, i.e.,
$y_1 (0) = y_0 (0) = O$ and
\begin{equation}
 \label{initial curve}
  \Gamma_D (0) 
  = L_1 (0) \cup L_0 (0)
  = \{ \lambda {\bf T}_0 ; \ \lambda \ge 0 \}.
\end{equation}
For evolution of a ``spiral'', 
a new facet should be generated as the resultant of
the evolution of present facets.
Let $T_1 = 0$ and inductively set
the generation time of facet $L_{k+1} (t)$ as
\[
 T_{k+1} = \sup \{ T > T_k ; \ d_k (t) < \rho_c \ell_k / U 
 \ \mbox{for} \ t \in [T_k, T] \}.
\]
When $t = T_{k+1}$, we add a new facet 
$L_{k+1} (T_{k+1})$
with $y_{k+1} (T_{k+1}) = O$ and $d_{k+1} (T_{k+1}) = 0$.
Then, change the spiral center to $y_{k+1} (t)$ from $y_k (t)$.

In summary, the algorithm of our discrete model for
evolving polygonal spiral by \eqref{geo mcf} is as follows:
\begin{enumerate}
 \renewcommand{\theenumi}{\Roman{enumi}}
 \item The generation time $T_k$ and 
       curve $\Gamma_D (T_k) = \bigcup_{j=0}^k L_j (T_k)$
       (with $d_k (T_k) = 0$) are given.
 \item Solve \eqref{ODE system1}--\eqref{ODE system2}
       on $[T_k, T_{k+1}]$
       to obtain the evolution of $\Gamma_D (t)$.
 \item When $t = T_{k+1}$, add a new facet 
       $L_{k+1} (T_{k+1})$
       with $y_{k+1} (T_{k+1}) = O$ 
       (then $d_{k+1} (T_{k+1}) = 0$)
       as the fixed center of $\Gamma_D (t)$.
       Then, return to (I).
\end{enumerate}
The existence and uniqueness 
of solution to \eqref{ODE system1}--\eqref{ODE system2},
the existence of the sequence $\{ T_k \}_{k=1}^\infty$
of the generation times,
$\lim_{k \to \infty} T_k = \infty$,
and the intersection-free result of $\Gamma_D (t)$
are obtained by \cite{IO:DCDS-B};
see it for details of the mathematical results.

\subsection{Level set method}
\label{sec: LV model}

We recall the level set method \cite{OTG:2015JSC}
for a evolving spiral corresponding to the discrete
model explained in the previous section.

Let $\Omega \subset \mathbb{R}^2$ be a bounded domain
with a smooth boundary.
Consider the evolution of a single spiral
by \eqref{geo mcf}, and 
have set the center of a spiral at the origin.
We give such a spiral curve and its direction 
of the evolution, which is denoted by ${\bf n} \in S^1$,
with the level set method due to \cite{OTG:2015JSC} as
\[
 \Gamma_L (t)
 = \{ x \in \overline{W} ; \ u(t,x) - \theta (x) \equiv 0
 \mod 2 \pi \mathbb{Z} \},
 \quad {\bf n} = - \frac{\nabla (u - \theta)}{|\nabla (u - \theta)|},
\]
where $W = \{ x \in \Omega ; \ |x| > \rho \}$
for a constant $\rho > 0$,
and $\theta = \arg x$.
According to \cite{Giga:2006}, 
we obtain the anisotropic curvature $H_\gamma$
of $\Gamma_L (t)$ as
\[
 H_\gamma = - \mathrm{div} \xi (- \nabla (u - \theta))
\]
with $\xi = D \gamma$ 
and 
$\gamma \in C^2 (\mathbb{R}^2 \setminus \{ 0 \})$
satisfying (\ref{assum: convexity})--(\ref{assum: positivity}).
It is well-known that 
\[
 \mathcal{W}_\gamma = \{ p \in \mathbb{R}^2 ; \ \gamma^\circ (p) \le 1 \}
\]
with $\gamma^\circ (p) = \sup \{ p \cdot q ; \ \gamma(q) \le 1 \}$,
and $H_\gamma = 1$ on $\mathcal{W}_\gamma$;
see \cite{Bellettini-Paolini:1996HMJ} for details.
Moreover, from the context of derivation of 
\eqref{ODE system1}--\eqref{ODE system2} as in \cite{IO:DCDS-B},
one can find a self-similar solution with extension of 
$\mathcal{W}_\gamma$ for the motion of closed curve by $V=1$
(\eqref{geo mcf} with $U=1$ and $\rho_c = 0$),
which means that we measure the normal velocity with 
the Finsler metric
\begin{align}
 & d_\gamma (x,y) = \gamma^\circ (x-y) \nonumber \\
 & \label{normalization}
 \mbox{with} \quad
 \gamma^\circ ({\bf N}_j) = 1 \quad \mbox{for} \ 
 j=0, 1, 2, \ldots, N_\gamma - 1.
\end{align}
Then, the normal velocity in this case should be given by
\[
 V_\gamma = \frac{u_t}{\gamma (- \nabla (u - \theta))}
\]
since $\gamma (D \gamma^\circ (p)) = 1$
for $p \in \mathbb{R}^2 \setminus \{ 0 \}$
under some additional regularity and convexity assumptions
on $\gamma$ and $\gamma^\circ$;
see \cite{Bellettini-Paolini:1996HMJ} for details.

As a boundary condition of the evolution with \eqref{geo mcf},
we impose the right angle condition between $\Gamma_L (t)$ and
$\partial W$.
Then, the level set equation
of the motion of spirals by \eqref{geo mcf}
is of the form
\begin{eqnarray}
 \label{lv mcf}
  \tilde{\beta} (\nabla (u - \theta))
  u_t - \tilde{\gamma} (\nabla (u - \theta))
  \left\{
   \rho_c \mathrm{div} \tilde{\xi} (\nabla (u - \theta))
   + U
  \right\}
  = 0 & \mbox{in} & (0,T) \times W, \\
 \label{lv nbc}
 \vec{\nu} \cdot \nabla (u - \theta) = 0
  & \mbox{on} & (0,T) \times W,
\end{eqnarray}
where $\vec{\nu} \in \mathbb{S}^1$ is the outer unit normal
vector field of $\partial W$,
and $\tilde{\beta} (p) = \beta (-p)$,
$\tilde{\gamma} (p) = \gamma (-p)$
and $\tilde{\xi} (p) = \xi (-p)$.
See \cite{Giga:2006} for details of the level set method.

Mathematical analysis for \eqref{lv mcf}--\eqref{lv nbc}
with $\gamma \in C^2 (\mathbb{R}^2 \setminus \{ 0 \})$ 
and $\beta \in C(\mathbb{R}^2 \setminus \{ 0 \})$
is established in \cite{Ohtsuka:2003wi}.
For given initial data $u_0 \in C(\overline{W})$,
there exists a unique global viscosity solution 
$u \in C([0,\infty) \times \overline{W})$
to \eqref{lv mcf}--\eqref{lv nbc} with $u(0,\cdot) = u_0$.
Moreover, the uniqueness of evolution of $\Gamma_L (t)$
is established in \cite{Goto:2008hy};
if there are continuous viscosity solutions 
$u$ and $v$ to \eqref{lv mcf}--\eqref{lv nbc}
satisfying $\Gamma_L^u (0) = \Gamma_L^v(0)$
with the same orientations, 
then $\Gamma_L^u (t) = \Gamma_L^v(t)$
for $t > 0$, where
$\Gamma_L^u (t) = \{ x \in \overline{W} ; \ 
u(t,x) - \theta (x) \equiv 0 \mod 2 \pi \mathbb{Z} \}$.
Hence, we may give an arbitrary $u_0 \in C(\overline{W})$
to obtain the motion of $\Gamma_L (t)$.
In this paper, we give $u_0$ for \eqref{initial curve}
as $u_0 \equiv \varphi_0$ due to \cite{OTG:2015JSC}.

Recall that we consider the situation such that
$\mathcal{W}_\gamma$ is a convex polygon.
The assumption (\ref{assum: crystalline}) 
is imposed for such a situation.
Then,
$\gamma$ is now given as
\begin{equation}
 \label{char form of gamma}
 \gamma (p) = \max_{0 \le j \le N_\gamma - 1} n_j \cdot p
 = \sum_{j=0}^{N_\gamma - 1} (n_j \cdot p) \chi_{Q_j} (p) 
\end{equation}
with some $Q_j \subset \mathbb{R}^2$ for $j=0,1,2, \ldots, N_\gamma - 1$,
where
\[
 \chi_Q (x) = \left\{
 \begin{array}{ll}
  1 & \mbox{if} \ x \in Q, \\
  0 & \mbox{otherwise}
 \end{array}
 \right.      
\]
for $Q \subset \mathbb{R}^2$.
The crucial problem for solving \eqref{lv mcf}--\eqref{lv nbc}
is how to treat $\mathrm{div} \tilde{\xi} (\nabla (u - \theta))$.
For this problem, approximation
of $\xi$ by the analogy of the stability result
as in \cite{Giga-Giga:2001} is a simple option.
From \eqref{char form of gamma}, we formally obtain
\[
 \xi (p) = \sum_{j=0}^{N_\gamma - 1} \chi_{Q_j} (p) n_j,
\]
so that we approximate $\chi_{Q_j}$
with the method as in \cite{ETT2005JCP}
to remove the singularities.
More precisely, we use the function
\[
 \sigma (z; p_1, p_2)
 := 
 \left\{ 
 \begin{array}{ll}
  {\displaystyle \frac{z}{\sqrt{z^2 + \varepsilon^2 (|p_1| + |p_2|)^2}} }
   & \mbox{if} \ z \neq 0, \\ [+5pt]
   0 & \mbox{otherwise}
 \end{array}
 \right.
\]
with $\varepsilon \ll 1$ to approximate the sign function $z/|z|$.
This function is also used in \cite{IO:DCDS-B}
when we approximate $\xi = D \gamma$ of 
$\gamma (p) = \| p \|_1 = |p_1| + |p_2|$
or $\gamma (p) = \| p \|_\infty = \max \{ |p_1|, |p_2| \}$
for $p = (p_1, p_2)$.
In general, 
consider the case 
when $Q$ is given as a level set of a continuous function $f$, i.e., 
$Q = \{x \in \mathbb{R}^2 ; \ f (x) > 0 \}$ and 
$\mathbb{R}^2 \setminus \overline{Q} 
= \{x \in \mathbb{R}^2 ; \ f (x) < 0 \}$.
Then, we approximate $\chi_Q$ by 
\[
 \chi_Q (x) \approx 
 \zeta (f(x); p_1, p_2), \quad
 \mbox{with} \ \zeta (z;p_1, p_2) := \frac{\sigma (z; p_1, p_2) + 1}{2} 
\]
for a suitable parameter $(p_1, p_2)$.
(We often choose $(p_1, p_2) = \nabla f$ 
like as in \cite{ETT2005JCP}, 
or $(p_1, p_2) = (1,0)$ for simplicity.)
Hence, we obtain the approximation 
\[
 \xi (p) \approx \sum_{j=0}^{N_\gamma - 1} 
 \zeta (f_j (x); p_1, p_2) n_j 
\]
by a level set functions $f_j \in C(\mathbb{R}^2)$ for $Q_j$.

\section{Measuring difference}

\subsection{Crystalline energy density}
\label{sec: how to define density}

Let us consider the situation such that
the Wulff shape $\mathcal{W}_\gamma$
and a support function $\gamma^\circ \colon \mathbb{R}^2 \to [0, \infty)$
satisfying
\[
 \mathcal{W}_\gamma
 = \{ p \in \mathbb{R}^2 ; \ \gamma^\circ (p) \le 1 \}
\]
are given.
Note that $\gamma^\circ = \sup \{ p \cdot q; \ \gamma (q) \le 1 \}$
is a convex and positively homogeneous of degree 1.
According to these facts and that $\mathcal{W}_\gamma$ is a convex polygon,
we assume that $\gamma^\circ$ is given as
\begin{equation}
 \label{support formula}
 \gamma^\circ (p) := 
 \max_{0 \le j \le N_\gamma - 1} m_j \cdot p, \quad
 m_j = \eta_j (\cos \psi_j, \sin \psi_j) 
\end{equation}
with $\eta_j > 0$ and $\psi_j \in \mathbb{R}$.
Assume that 
\begin{enumerate}
 \renewcommand{\theenumi}{$\gamma$\arabic{enumi}}
 \item $\psi_0 < \psi_1 < \psi_2 \cdots < \psi_{N_\gamma - 1} 
       < \psi_0 + 2\pi$,
 \item $\psi_j < \psi_{j+1} < \psi_j + \pi$
       for $j = 0, 1, 2, \ldots, N_\gamma-1$.
 \item \label{assum: partition}
       $P_j = \{ p \in \mathbb{R}^2 ; \ 
       m_j \cdot p \ge m_k \cdot p \ 
       \mbox{for} \ k = 0, 1, 2, \ldots, N_\gamma - 1 \}
       = \Xi_{j,j-1} \cap \Xi_{j,j+1} \neq \emptyset$
       for $j=0,1,2,\ldots, N_\gamma - 1$,
       where $\Xi_{j,k}
       = \{ p \in \mathbb{R}^2 ; \ 
       m_j \cdot p \ge m_k \cdot p \}$.
\end{enumerate}
(Note that $\psi_{j + nN_\gamma} = \psi_j$
for $n \in \mathbb{Z}$.)
We now propose a practical way to reconstruct
a convex and piecewise linear
$\gamma \colon [0,\infty) \to [0,\infty)$
from the above settings.
Note that we do not impose the normalizing assumption
\eqref{normalization} in this section.

We first remark that
\[
 \mathcal{F}_\gamma
 = \{ p \in \mathbb{R}^2 ; \ \gamma (p) \le 1 \}
 = \{ p \in \mathbb{R}^2 ; \ p \cdot q \le \gamma^\circ (q) 
 \quad \mbox{for} \ q \in \mathbb{S}^1 \}
\]
when $\gamma$ is convex,
since $\gamma = (\gamma^\circ)^\circ$.
Then, by (\ref{assum: partition}), we find $\theta_0, \theta_1, \theta_2,
\ldots, \theta_{N_\gamma - 1}$ such that
$\theta_j < \theta_{j+1} < \theta_j + 2 \pi$ and
\[
 \gamma^\circ (q) = m_j \cdot q
 \quad \mbox{if} \ q = (\cos \theta, \sin \theta) \
 \mbox{with} \ \theta \in [\theta_j, \theta_{j+1}]
\]
for $j=0,1,2,\ldots, N_\gamma - 1$.
It should be calculated by
\[
 m_j \cdot (\cos \theta_j, \sin \theta_j)
 = m_{j-1} \cdot (\cos \theta_j, \sin \theta_j).
\]
In fact, by \eqref{support formula} we have
\begin{equation}
 \label{formula: support 1}
 \cos \theta_j (\eta_j \cos \psi_j - \eta_{j-1} \cos \psi_{j - 1})
 + \sin \theta_j (\eta_j \sin \psi_j - \eta_{j-1} \sin \psi_{j-1}) = 0. 
\end{equation}
Let $a_j$, $b_j$ be constants defined by 
\begin{align*}
 & a_j = \eta_j \cos \psi_j - \eta_{j-1} \cos \psi_{j - 1}, \quad 
 b_j = \eta_j \sin \psi_j - \eta_{j-1} \sin \psi_{j-1}, 
\end{align*}
and $c_j$ be a constant satisfying
\[
 \cos c_j = \frac{a_j}{\sqrt{a_j^2 + b_j^2}}, \quad
 \sin c_j = \frac{b_j}{\sqrt{a_j^2 + b_j^2}}.
\]
Then, \eqref{formula: support 1} yields that
\begin{equation}
 \label{def of theta_j}
 \cos (\theta_j - c_j) = 0, \quad
 \mbox{i.e.}, \quad
 \theta_j = c_j + \frac{\pi}{2}. 
\end{equation}

Let us consider the formula
$p \cdot q \le \gamma^\circ (q)$
with $q = (\cos \theta, \sin \theta)$
and $p = (x,y)$.
If $\theta \in [\theta_j, \theta_{j+1}]$,
then we observe that
\[
 \cos \theta (x - \eta_j \cos \psi_j) + \sin \theta 
 (y - \eta_j \sin \psi_j) \le 0
 \quad \mbox{for} \ \theta \in [\theta_j, \theta_{j+1}].
\]
Then, one can find that
\begin{align*}
 & \{ (x,y) \in \mathbb{R}^2 ; \ 
 \cos \theta (x - \eta_j \cos \psi_j) + \sin \theta 
 (y - \eta_j \sin \psi_j) \le 0
 \quad \mbox{for} \ \theta \in [\theta_j, \theta_{j+1}]
 \} \\
 & =
 \Pi_{j,j} \cap \Pi_{j,j+1}, 
\end{align*}
where
\[
 \Pi_{j,k} 
 = \{ (x,y) \in \mathbb{R}^2 ; \ 
 \cos \theta_k (x - \eta_j \cos \psi_j) + \sin \theta_k
 (y - \eta_j \sin \psi_j) \le 0\}.
\]
Moreover, 
one can find $\Pi_{j,j+1} = \Pi_{j+1, j+1}$.
In fact, by definition of $\Pi_{j, j+1}$ and 
\eqref{def of theta_j}, we observe that
\begin{align*}
 & \cos \theta_{j+1} (x - \eta_j \cos \psi_j) + \sin \theta_{j+1}
 (y - \eta_j \sin \psi_j) \\
 & = 
 \cos \theta_{j+1} (x - \eta_{j+1} \cos \psi_{j+1}) + \sin \theta_{j+1}
 (y - \eta_{j+1} \sin \psi_{j+1}) \\
 & \qquad 
 + \cos \theta_{j+1} (\eta_{j+1} \cos \psi_{j+1}- \eta_j \cos \psi_j) 
 + \sin \theta_{j+1} (\eta_{j+1} \sin \psi_{j+1} - \eta_j \sin \psi_j)
 \\
 & = 
 \cos \theta_{j+1} (x - \eta_{j+1} \cos \psi_{j+1}) + \sin \theta_{j+1}
 (y - \eta_{j+1} \sin \psi_{j+1}) \\
 & \qquad 
 + a_{j+1} \cos \theta_{j+1} + b_{j+1} \sin \theta_{j+1} \\
 & = 
 \cos \theta_{j+1} (x - \eta_{j+1} \cos \psi_{j+1}) + \sin \theta_{j+1}
 (y - \eta_{j+1} \sin \psi_{j+1}) \\
 & \qquad + \sqrt{a_{j+1}^2 + b_{j+1}^2} 
 \cos (\theta_{j+1} - c_{j+1}) \\
 & = 
 \cos \theta_{j+1} (x - \eta_{j+1} \cos \psi_{j+1}) + \sin \theta_{j+1}
 (y - \eta_{j+1} \sin \psi_{j+1}),
\end{align*}
which implies $\Pi_{j,j+1} = \Pi_{j+1, j+1}$.
Hence, we obtain
\begin{align*}
 \mathcal{F}_\gamma
 = \{ p \in \mathbb{R}^2 ; \ 
 p \cdot q \le \gamma^\circ (q) \ \mbox{for} \ q \in \mathbb{S}^1 \}
 = \bigcap_{j=0}^{N_\gamma - 1} \Pi_{j,j}.
\end{align*}
Set
$r_j = \left[\eta_j 
(\cos \theta_j \cos \psi_j + \sin \theta_j \sin \psi_j) \right]^{-1}
= \left[ \eta_j \cos (\theta_j - \psi_j) \right]^{-1}$,
and 
\[
 \gamma (p) = \max_{0 \le j \le N_\gamma -1} n_j \cdot p, \quad
 n_j = r_j (\cos \theta_j, \sin \theta_j).
\]
Then, we observe that
\[
 \cos \theta_j (x - \eta_j \cos \psi_j) + \sin \theta_j
 (y - \eta_j \sin \psi_j)
 = \frac{n_j \cdot p}{r_j} - \frac{1}{r_j}
 \quad \mbox{with} \ p = (x,y),
\]
which 
implies that
\begin{align*}
 \bigcap_{j=0}^{N_\gamma - 1} \Pi_{j,j}
 = \bigcap_{j=0}^{N_\gamma - 1} 
 \{ p \in \mathbb{R}^2 ; \ n_j \cdot p \le 1 \}
 = 
 \{ p \in \mathbb{R}^2 ; \ \max_{0 \le j \le N_\gamma -1} n_j \cdot p \le 1 \}.
\end{align*}
Hence, we observe that
$\gamma (p) = \max_{0 \le j \le N_\gamma - 1} n_j \cdot p$.

\par
\bigskip
\noindent
\textbf{Summary.}
Assume that ($\gamma$1)--($\gamma$3) hold.
Let $\gamma^\circ \colon \mathbb{R}^2 \to [0,\infty)$ 
be given as
\[
 \gamma^\circ (p) = \max_{0 \le j \le N_\gamma - 1} m_j \cdot p
 \quad \mbox{with} \ m_j = \eta_j (\cos \psi_j, \sin \psi_j).
\]
Set 
\begin{itemize}
 \item $\theta_j = c_j + \pi / 2$ with $c_j \in \mathbb{R}$
       such that
       \begin{align*}
	& \cos c_j = \frac{a_j}{\sqrt{a_j^2 + b_j^2}}, \quad
	\sin c_j = \frac{b_j}{\sqrt{a_j^2 + b_j^2}}, \\
	& a_j = \eta_j \cos \psi_j - \eta_{j-1} \cos \psi_{j-1}, \quad
	b_j = \eta_j \sin \psi_j - \eta_{j-1} \sin \psi_{j-1},
       \end{align*}
 \item $r_j = \left[ \eta_j \cos (\theta_j - \psi_j) \right]^{-1}$
\end{itemize}
for $j=0, 1, 2, \ldots, N_\gamma - 1$.
Then,
\[
 \gamma (p) = \max_{0 \le j \le N_\gamma - 1} n_j \cdot p
 \quad \mbox{with} \ n_j = r_j (\cos \theta_j, \sin \theta_j).
\]

\par
\bigskip
\noindent
\begin{remark}
 \begin{enumerate}  
  \item When we give only the parameters of $\ell_j$ and $\psi_j$
	for $\mathcal{W}_\gamma$,
	then we have to set the location
	of the origin $O \in \mathcal{W}_\gamma$
	to determine $\gamma^\circ$.
	Note that 
	$\eta_j$ depends on the location of the origin 
	in $\mathcal{W}_\gamma$.
  \item 
 There is a case that
 $P_j = \emptyset$ and thus
 $P_k \neq \Xi_{k,k-1} \cap \Xi_{k,k+1}$
 for some $j,k \in \{ 0,1,2, \ldots, N_\gamma - 1 \}$ 
 when $N_\gamma \ge 4$,
 even if $\psi_j$ satisfies ($\gamma$1).
 In fact, $\gamma^\circ (p) = \max_{0 \le j \le 3} n_j \cdot p$ with
 \[
 n_0 = (3,0), \ n_1 = (1,1), \ n_2 = (0,2), \ n_3 = (-1,-1)
 \]
 implies that $P_1 = \emptyset$.  
  \item Notice that the above way also can be applied to
	construct $\gamma^\circ$ from a given $\gamma$.
 \end{enumerate}
\end{remark}

\subsection{Difference function}
\label{sec: diff function}

\begin{figure}[htbp]
 \begin{center}
  \includegraphics[scale=0.2]{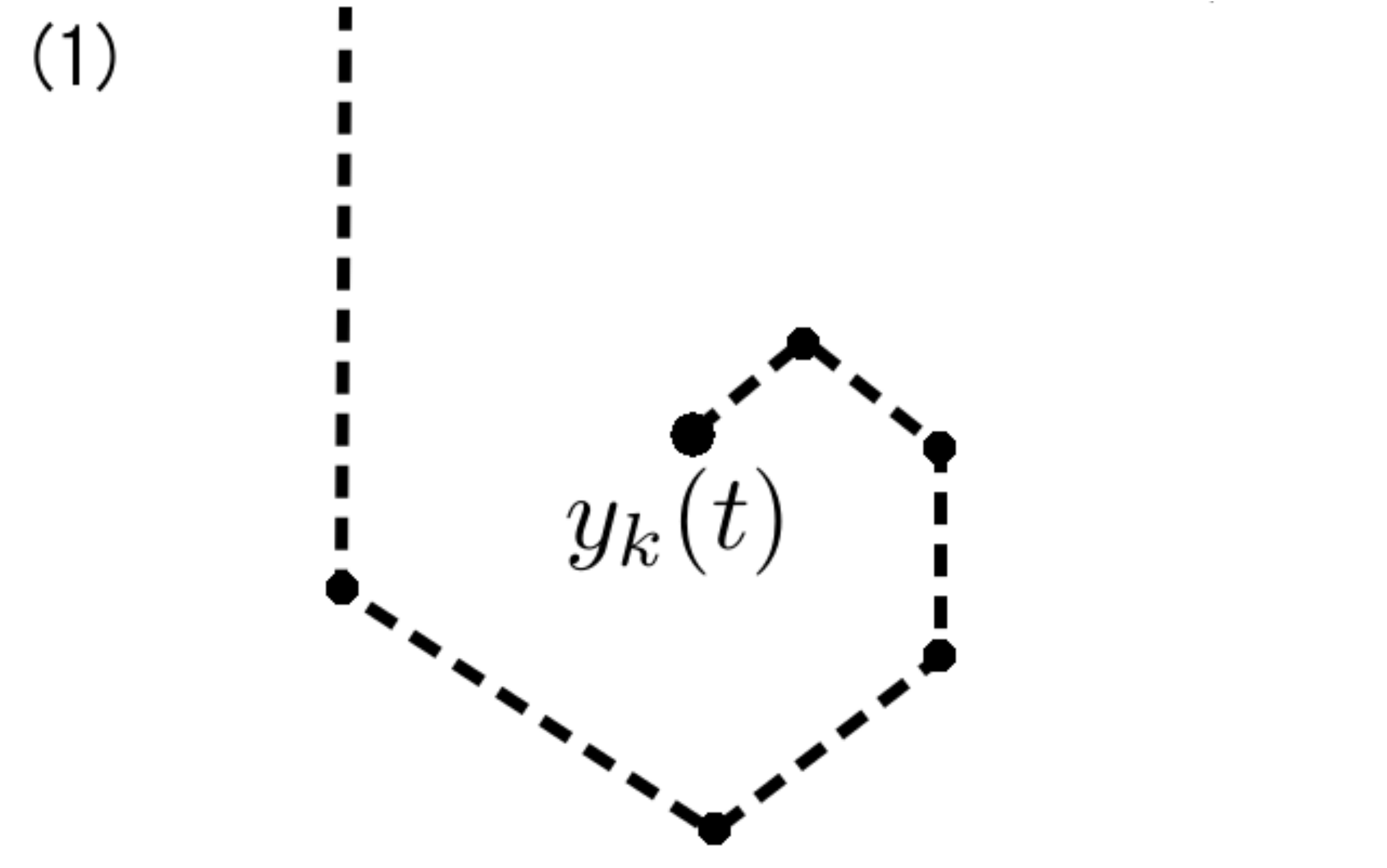}
  \includegraphics[scale=0.2]{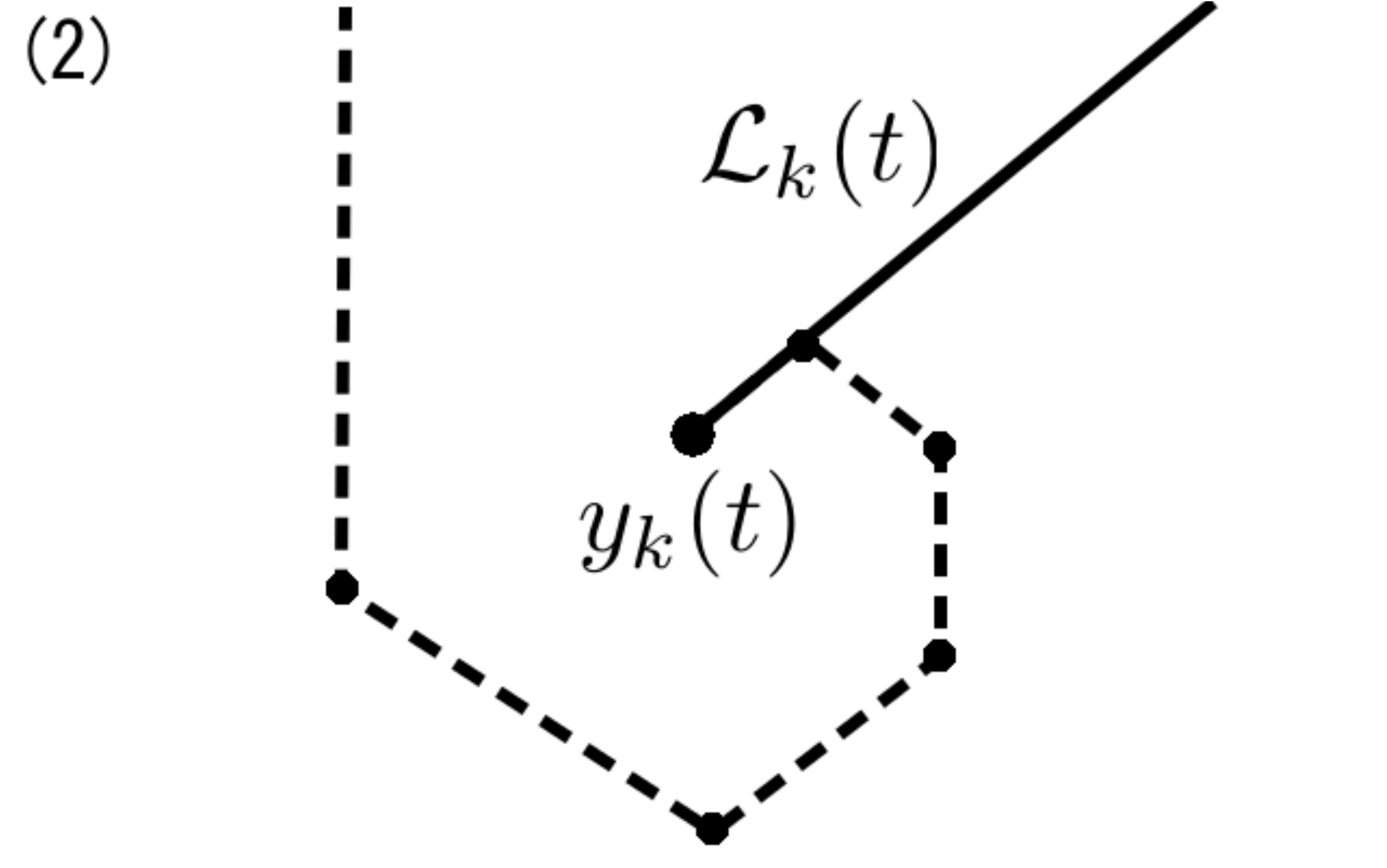} \\
  \includegraphics[scale=0.2]{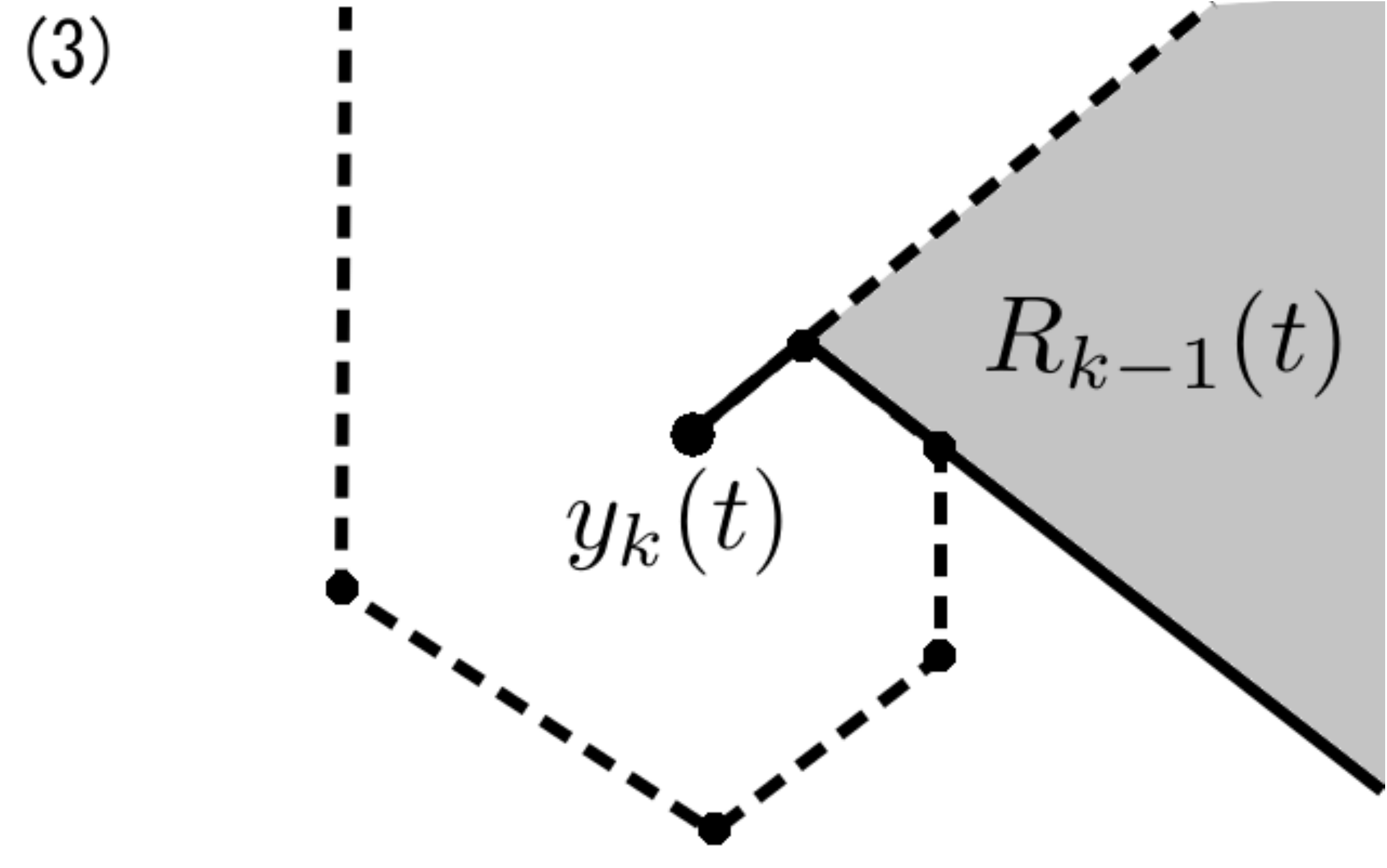}
  \includegraphics[scale=0.2]{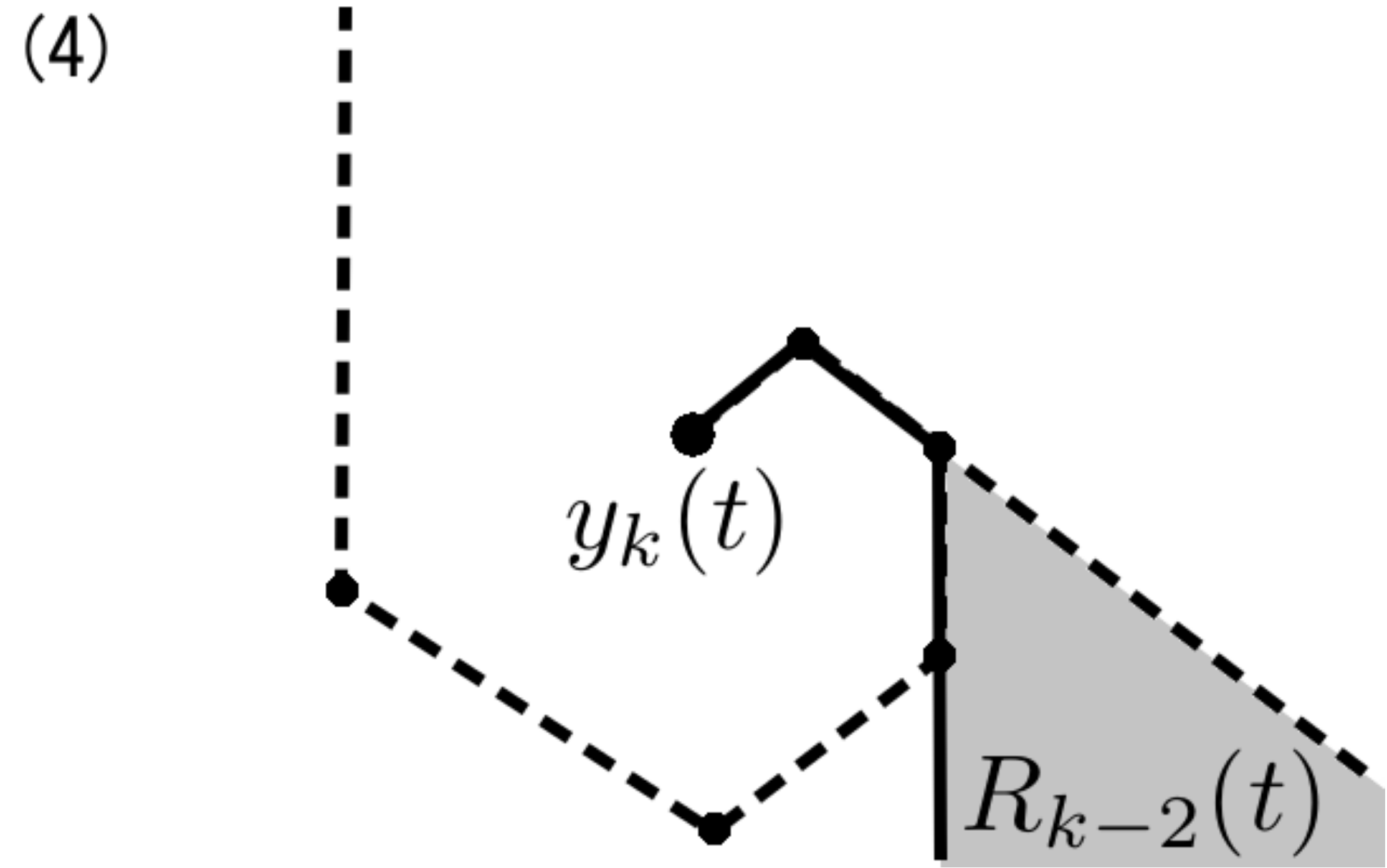}
  \caption{Construction of $\theta_D (t, x)$;
  we construct a branch of $\arg x$ 
  whose discontinuities
  are only on $\Gamma (t)$(dashed line in (1)).
  For this purpose we first construct
  $\vartheta (x) = \arg x$ whose discontinuities
  are only on $\mathcal{L}_k (t)$
  (solid line in (2)).
  Then, we make go down 
  the height of $\vartheta (x)$
  on $R_{j} (t)$ (gray region in (3) or (4))
  with the jump-height $2 \pi$
  from $j=k-1$ to $j=0$ inductively
  to remove illegal discontinuities.
  The solid line in figure (3) or (4)
  denotes the discontinuity
  of $\Theta_{k,k-1}$ or $\Theta_{k,k-2}$,
  respectively.}
  \label{height function for line segments}
 \end{center}
\end{figure}

Once we obtain
$\Gamma_D (t) = \bigcup_{j=0}^k L_j (t)$
or
$\Gamma_L (t)$,
then we compare with $\Gamma_D (t)$ and $\Gamma_L (t)$
by calculating the measure of the region 
enclosed by $\Gamma_D (t)$ and $\Gamma_L (t)$.
It is established as follows;
we construct the height functions
\[
 h_D (t, x) = \frac{1}{2 \pi} \theta_D (t,x), \quad
 h_L (t, x) = \frac{1}{2 \pi} \theta_L (t,x)
\]
of the stepwise surface 
at $\Gamma_D (t)$ or $\Gamma_L (t)$
with step height $h_0 = 1$, respectively.
Note that $\theta_D (t,x)$ or $\theta_L (t,x)$
is a branch of $\theta (x)$
whose discontinuity is only on $\Gamma_D (t)$ or $\Gamma_L (t)$,
respectively.
According to \cite{OTG:2015JSC}, 
the practical way to construct $\theta_L (t,x)$
is given in \cite{OTG:2015JSC}.
Hence, we here give a practical way to construct
$\theta_D (t,x)$.
\begin{enumerate}
 \item We first pick up the rotation number
       $n \in \mathbb{N}$
       for the facet number $k \in \mathbb{N}$
       of $\Gamma (t) = \bigcup_{j=0}^k L_j (t)$,
       i.e., 
       $k = \bar{k} + n N_\gamma$ with
       $\bar{k} \in \{ 0, 1, 2, \ldots, N_\gamma - 1 \}$.
 \item Then, we now set
       \[
	\Theta_k (x) = \arg (x)
       \in [\varphi_{\bar{k}} + 2 \pi n - \pi / 2,
       \varphi_{\bar{k}} + 2 \pi (n+1) - \pi / 2);
       \]
       a branch of $\arg x$ whose discontinuity is
       only on
       \[
       \mathcal{L}_k (t)
       = \{ r {\bf T}_k ; \ r > 0 \}.
       \]
       (See Figure \ref{height function for line segments}(2).)
 \item Let us set
       \[
       R_{k-1} (t) := \{ x \in \mathbb{R}^2 ; \
       x \cdot {\bf N}_k < s_k (t), \
       x \cdot {\bf N}_{k-1} \ge s_{k-1} (t) \}
       \]
       (gray regions in Figure
       \ref{height function for line segments}(3)).
       To remove a discontinuity on a dash line
       in $\partial R_{k-1} (t)$,
       we set
       \[
       \Theta_{k,k-1} (x)
       = \Theta_k (x) - 2 \pi \chi_{R_{k-1} (t)} (x).
       \]
 \item We inductively set
       \begin{align*}
	\Theta_{k,k-\ell} (x)
	& = \Theta_{k,k-\ell+1} (x)
	- 2 \pi \chi_{R_{k - \ell} (t)} (x) \\
	& = \Theta_k (x) - 2 \pi \sum_{j=1}^\ell
	\chi_{R_{k-j} (t)} (x)	
       \end{align*}
       to remove illegal discontinuities
       of $\Theta_{\ell-1}$
       from $\ell = 1$ to $\ell = k$,
       where
       \[
	R_j (t) := \{ x \in \mathbb{R}^2 ; \
       x \cdot {\bf N}_{j+1} < s_{j+1} (t), \
       x \cdot {\bf N}_j \ge s_j (t) \}
       \]
       for $j=0,1,\ldots, k-1$
       (see Figure \ref{height function for line segments}(4)
       for $R_{k-2} (t)$).
\end{enumerate}
Consequently, we set
\[
 \theta_D (t,x)
 = \Theta_{k,0} (x)
 = \Theta_k (x)
 - 2 \pi \sum_{j=0}^{k-1} \chi_{R_j (t)} (x),
 \quad h_D (t,x) = \frac{1}{2\pi} \theta_D (t,x).
\]
Hence, we can define
the difference $\mathcal{D} (t)$
between $\Gamma_D (t)$ and $\Gamma_L (t)$
as
\begin{equation}
 \label{area difference}
 \mathcal{D} (t)
 = \frac{1}{|W|}
 \int_{W}
 |h_D (t, x) - h_L (t,x)| \textup{d}x.
\end{equation}

\section{Numerical results}
\label{sec: numerical results}

In this section, we present some
numerical simulations measuring 
the difference between $\Gamma_D (t)$
and $\Gamma_L (t)$
evolving by
\begin{align*}
 V_\gamma = 1 - \rho_c H_\gamma, 
\end{align*}
i.e., \eqref{geo mcf}
with $\beta \equiv 1$ and $U = 1$
for some kinds of $\gamma$.
The initial curve is chosen as
\[
 \Gamma_D (0) = \Gamma_L (0) = L_0 (0)
 = \{ \lambda {\bf T}_0 ; \ \lambda > 0 \},
\]
and then \eqref{initial curve} for the discrete model,
and $u (0,x) = \varphi_0$ for the level set method.
Throughout this section, we set
\[
 \Omega = [-1.5,1.5]^2, \
 D_s = \{ x_{i,j} = (i \Delta x, j \Delta x) ; \
 -75s \le i,j \le 75s \}
\]
for some $s \in \mathbb{N}$,
and then $\Delta x = 0.02/s$.
In the following subsections, we will choose time intervals of
the numerical simulations so that the curves $\Gamma_D (t)$
does not touch to the outer boundary $\partial \Omega$.
In other words, we avoid the situation that
the boundary condition on $\partial \Omega$
makes difference between $\Gamma_D (t)$ and $\Gamma_L (t)$.
Note that, however, the difference of the boundary condition
at the center and the evolution law of the first facet $L_0 (t)$
are still remains.

We calculate the ODE system \eqref{ODE system1}--\eqref{ODE system2}
by 4-th order Runge-Kutta method with the time span 
$\Delta t = 10^{-6}$.
From these numerical results, we construct $h_D (t,x)$
on each numerical mesh $D_s$ to compare the results
with those from the level set method.
On the other hand, the level set equation
\eqref{lv mcf}--\eqref{lv nbc} is calculated by
the explicit finite difference scheme as in \cite{OTG:2015JSC}
with the time span
$\Delta t = 0.1 \times \Delta x^2$.
See also \cite{OTG:2015JSC} for the way to construct
$h_L (t,x)$ with the step height $h_0 = 1$.
To draw a graph of $\mathcal{D} (t)$, we pick up the data
$\mathcal{D} (t_k) = \mathcal{D} (kT/20)$ for $0 \le k \le 20$
on the calculating time interval $[0,T]$.

We now recall the difference between the discrete model
in \S \ref{sec: ODE model} 
and the level set method in \S \ref{sec: LV model}.
\begin{enumerate}
 \item The domain of the level set method has a
       ``center'' $B_\rho = \{ x \in \mathbb{R}^2 ; \ |x| \le \rho \}$
       with a finite radius $\rho > 0$.
       However, the discrete has the center at the origin
       as a point (null set).
 \item The boundary conditions are different:
       \begin{itemize}
	\item {[Discrete model]}
	      $L_0 (t)$ evolves by $V = 1$ since $d_0 (t) = \infty$.
	      On the other hand, the behavior of the facets 
	      associated with center is imposed with
	      fixing and the generation rule of new facets.
	\item {[Level set method]}
	      The right angle conditions,
	      in particular, $\Gamma_L(t) \perp \partial B_{\rho}$
	      and $\Gamma_L (t) \perp \partial \Omega$
	      are imposed by \eqref{lv nbc}.
       \end{itemize}
\end{enumerate}
Because of the above differences, 
we have no conjectures of convergence between
$\Gamma_L (t)$ and $\Gamma_D (t)$ now.
Moreover, from the numerical results of the isotropic case
in \cite{OTG:2015JSC, OTG:2018CGD}, not only tending 
the approximation parameters to zero but also
letting $\rho \to 0$ is required for
numerical accuracy.
Thus, we shall check the numerical results with
fixed radius $0 < \rho \ll 1$ and 
reducing radius $\rho = O(\Delta x)$.

\subsection{Square spiral}
\label{sec: square}

\begin{figure}[htbp]
 \begin{center}
  \includegraphics[scale=0.55]{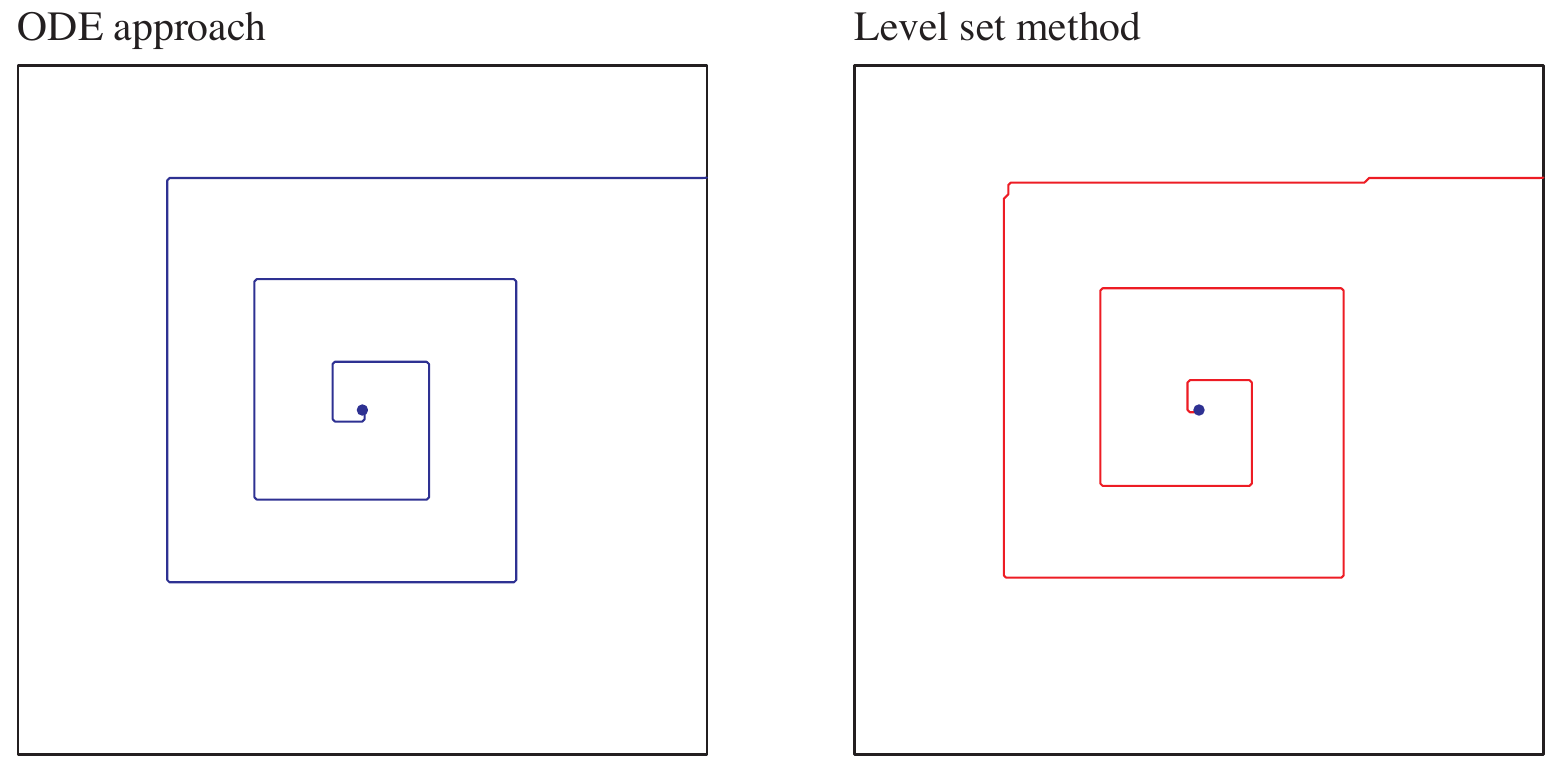}
  \caption{Profiles of the square spiral at $t=1$.
  The level set method is calculated with
  $\rho = 0.02 - 10^{-8}$ and $\Delta x = 0.0050$.}
  \label{profile cube2}
 \end{center}
\end{figure}

The first examination is the square spiral case, i.e., 
\[
 \mathcal{W}_\gamma = \{ p = (p_1, p_2) ; \ 
 \max \{ |p_1|, |p_2| \} \le 1 \}.
\]
Thus, we define the parameters of $\mathcal{W}_\gamma$
for the discrete model as
\[
 \varphi_j = \frac{\pi j}{2}, \quad
 \ell_j = 2
 \quad \mbox{for} \ j=0,1,2,3.
\]
For the level set equation, since 
$\gamma^\circ (p) = \max \{ |p_1|, |p_2| \}$
for $p = (p_1, p_2)$, we observe that
\[
 \gamma (p) = |p_1| + |p_2|, \quad
 \mbox{then} \ 
 \xi (p) = (\mathrm{sgn} (p_1), \mathrm{sgn} (p_2)).
\]
We calculate the ODE system \eqref{ODE system1}--\eqref{ODE system2}
and the level set equation \eqref{lv mcf}--\eqref{lv nbc}
for 
\begin{equation}
 \label{geo eq: cube-diag}
  V = 1 - 0.02 H_\gamma
  \quad (\mbox{i.e.,} \ \beta \equiv 1, \ U=1 \ \mbox{and} \ \rho_c = 0.02)
\end{equation}
on the time interval $[0,1]$.
See \cite[\S 4]{IO:DCDS-B} for details to approximate 
$\xi$ of the above $\gamma$.
Figure \ref{profile cube2} are profiles of the diagonal spiral
at $t=1$ with the above setting.
Note that, in this and following sections,
the profile of spirals by the level set method 
is calculated with
$\rho = 0.02 - 10^{-8}$ and $\Delta x = 0.0050$ ($s=4$).

\begin{figure}[htbp]
 \begin{center}
  \includegraphics[scale=0.5]{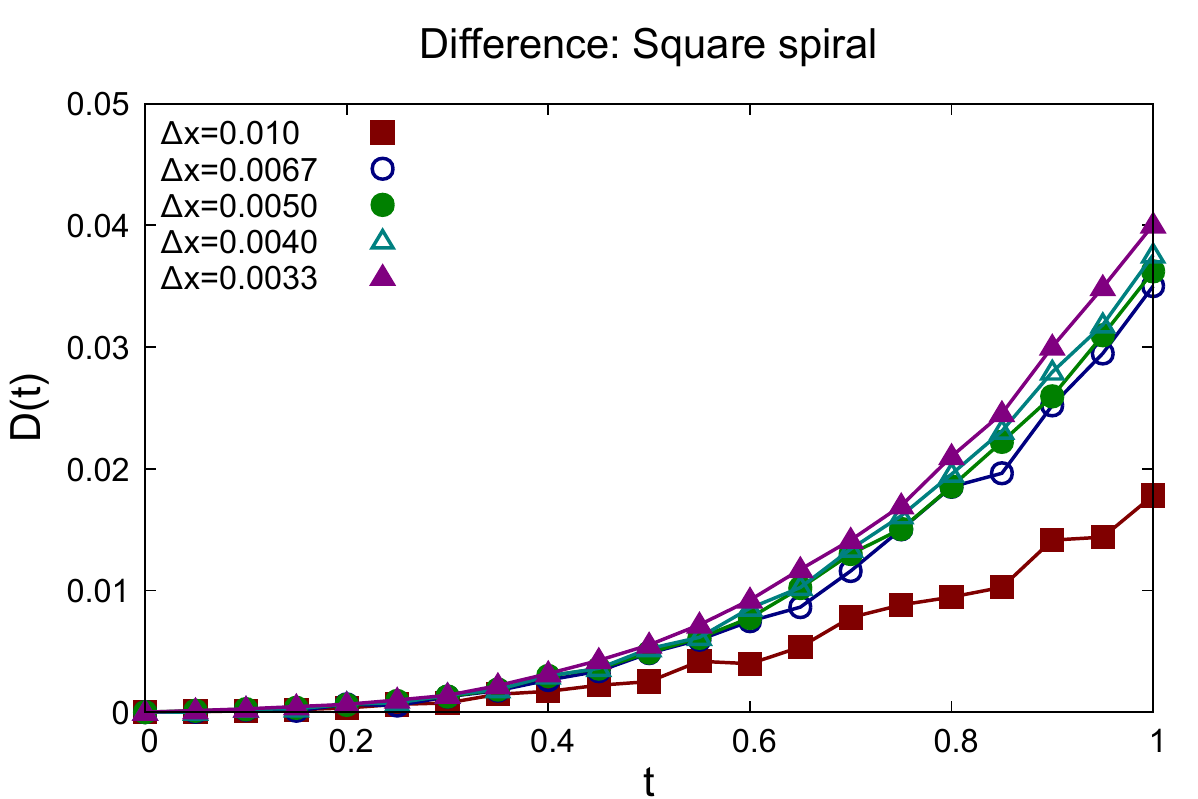}
  \includegraphics[scale=0.5]{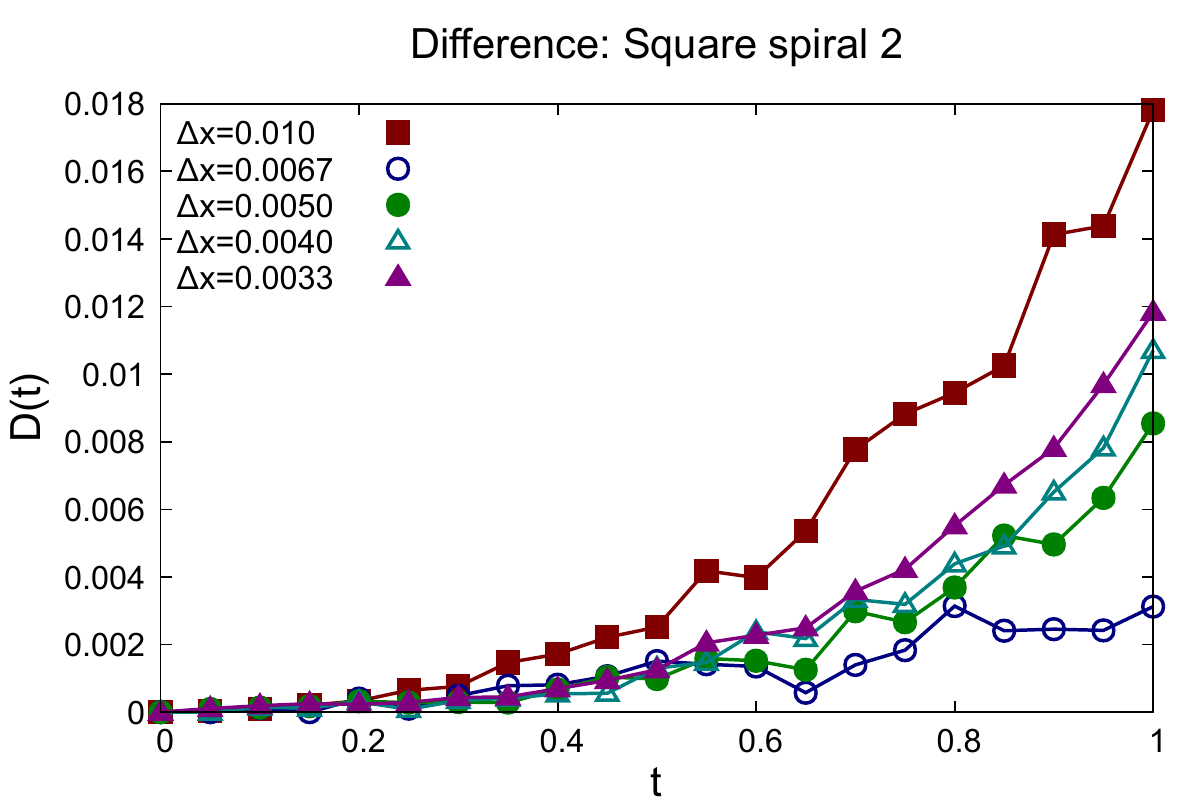}
  \caption{Graphs of functions $t \mapsto \mathcal{D} (t)$
  for the square spiral with a fixed center radius
  $\rho = 0.02 - 10^{-8}$(left), 
  and with a reduced center radius
  $\rho = (2 - 10^{-8}) \Delta x$(right).
  }
  \label{diff square}
 \end{center}
\end{figure}

The left figure of
Figure \ref{diff square} presents the graph of $\mathcal{D} (t)$
for $s=2,3,4,5,6$ with a fixed center radius
$\rho = 0.02 - 10^{-8}$.
One can find that the differences are
less than $4 \%$ of the area $|W|$
for all cases,
although the value of $\mathcal{D} (t)$
becomes worse when we choose smaller $\Delta x$.
The best one is the case with $\Delta x = 0.010$ ($s=2$).

On the other hand, we obtain better results when
$\rho = O(\Delta x)$.
The right figure of
Figure \ref{diff square} presents the graph of $\mathcal{D} (t)$
for $s=2,3,4,5,6$ with the center size 
$\rho = (2 - 10^{-8}) \Delta x$, i.e., 
the setting $\rho \to 0$ as $\Delta x \to 0$.
Note that the cases of $\Delta x = 0.010$ ($s=2$) in both figures of
Figure \ref{diff square}
is the same.
One can find that the differences are less than 
2.5\% of the area $|W|$ for all cases, and 
$\mathcal{D}(t)$ of the cases with $s \ge 3$ are
smaller than that of $s=2$,
although the smallest $\mathcal{D} (t)$
is the case $\Delta x = 0.0067$ ($s=3$).

\subsection{Diagonal spiral}
\label{sec: diagonal}

\begin{figure}[htbp]
 \begin{center}
  \includegraphics[scale=0.55]{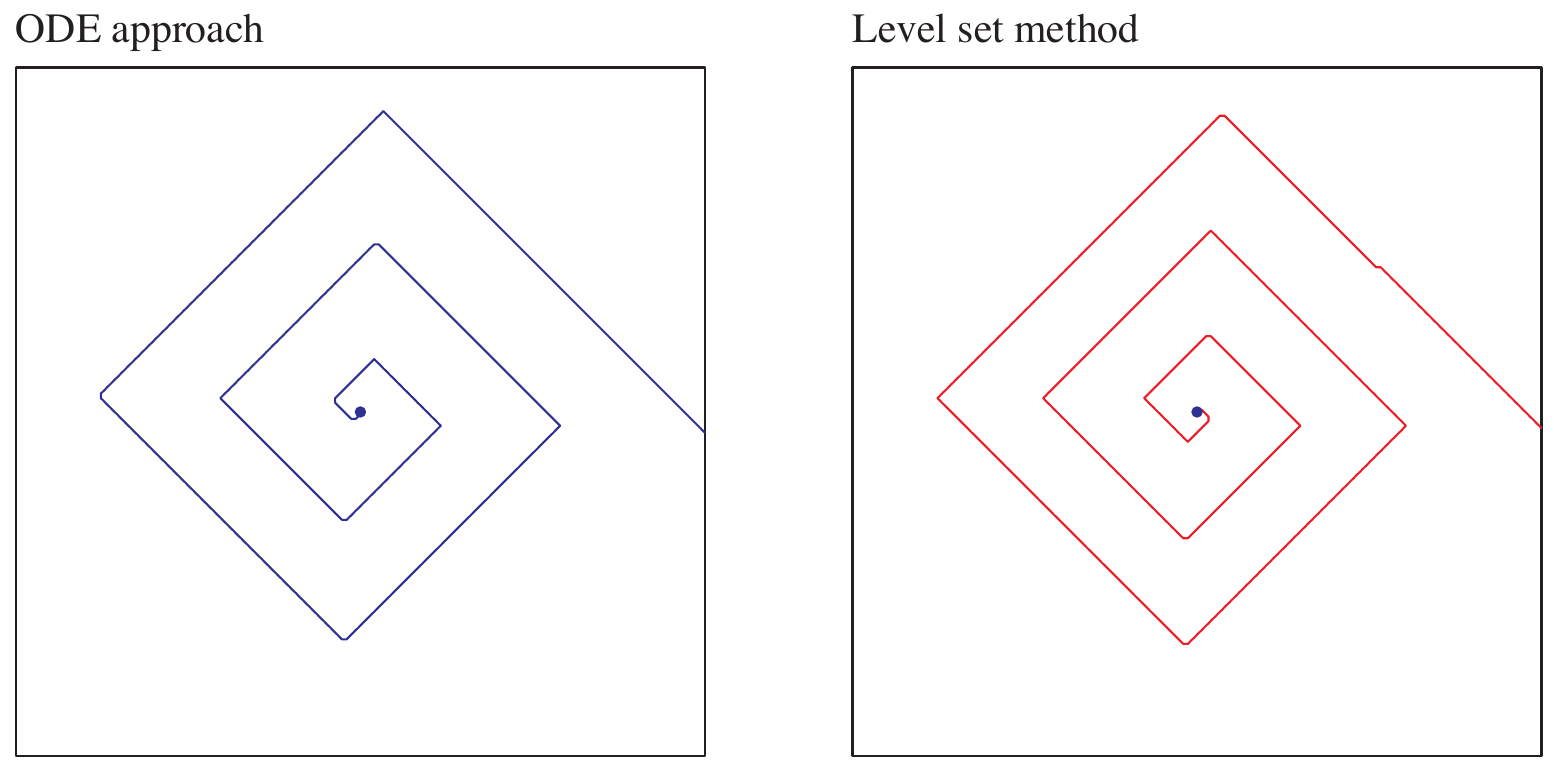}
  \caption{Profiles of the diagonal spiral at $t=1$.
  The level set method is calculated with
  $\rho = 0.02 - 10^{-8}$ and $\Delta x = 0.0050$.}
  \label{profile diag2}
 \end{center}
\end{figure}

The second examination is the diagonal spiral case, i.e.,
the $\pi/4$ rotation of the first case;
and then
\[
 \varphi_j = \frac{\pi j}{2} + \frac{\pi}{4}, \quad
 \ell_j = 2
 \quad \mbox{for} \ j=0,1,2,3
\]
for the discrete model.
In this case, one can find that
\[
 \mathcal{W}_\gamma = \{ p = (p_1, p_2) ; \ 
 |p_1|+|p_2| \le \sqrt{2} \},
\]
and thus
\[
 \gamma^\circ (p) = \frac{|p_1| + |p_2|}{\sqrt{2}}.
\]
For the level set equation, 
we set
\[
 \gamma (p) = \sqrt{2} \max \{ |p_1|, |p_2| \}.
\]
According to \cite{IO:DCDS-B}, it is represented as
\[
 \gamma (p)
 = \left| \frac{p_1 + p_2}{\sqrt{2}} \right|
 + \left| \frac{p_1 - p_2}{\sqrt{2}} \right|
\]
and thus
\[
 \xi (p)
 = \frac{1}{\sqrt{2}}
 (\mathrm{sgn} (p_1 + p_2) + \mathrm{sgn} (p_1 - p_2),
 \mathrm{sgn} (p_1 + p_2) - \mathrm{sgn} (p_1 - p_2)).
\]
See \cite{IO:DCDS-B} for the approximation of the above $\xi$.
We calculate the ODE system \eqref{ODE system1}--\eqref{ODE system2}
and the level set equation \eqref{lv mcf}--\eqref{lv nbc}
for \eqref{geo eq: cube-diag} on the time interval $[0,1]$.
Figure \ref{profile diag2} are profiles of the diagonal spiral
at $t=1$ with the above setting.

\begin{figure}[htbp]
 \begin{center}
  \includegraphics[scale=0.5]{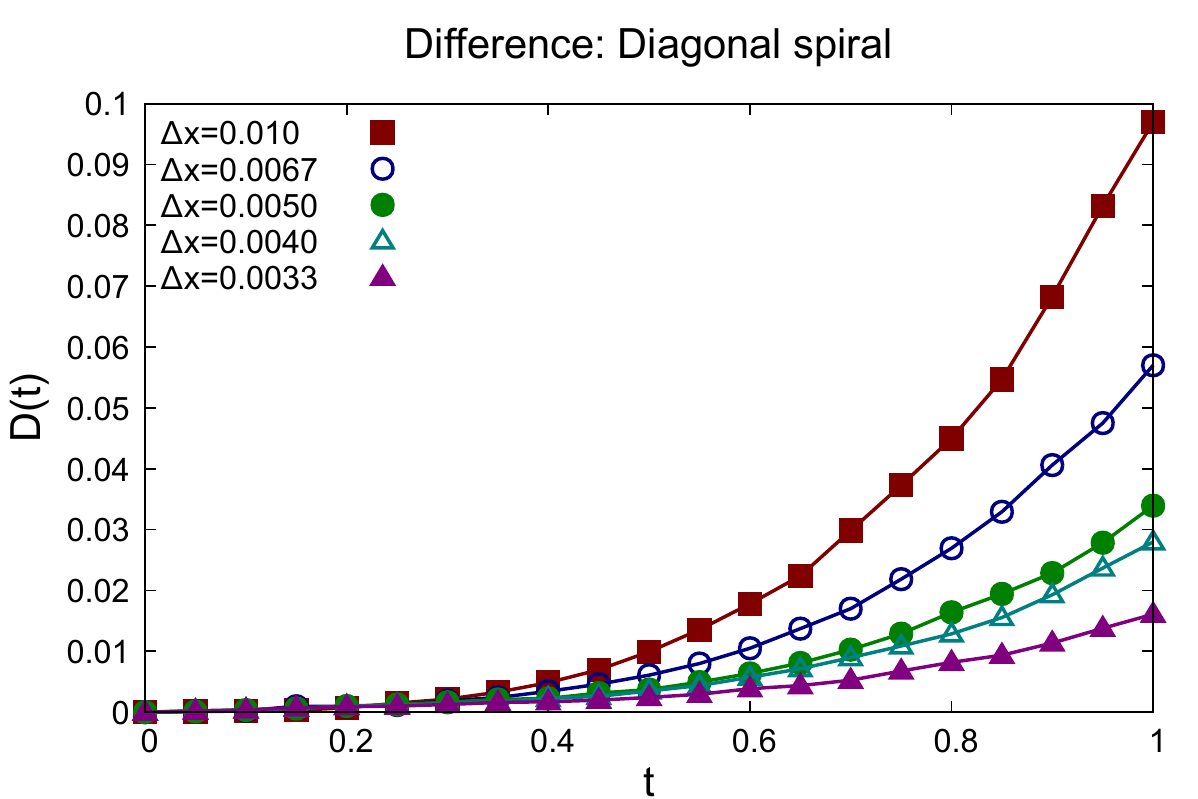}
  \includegraphics[scale=0.5]{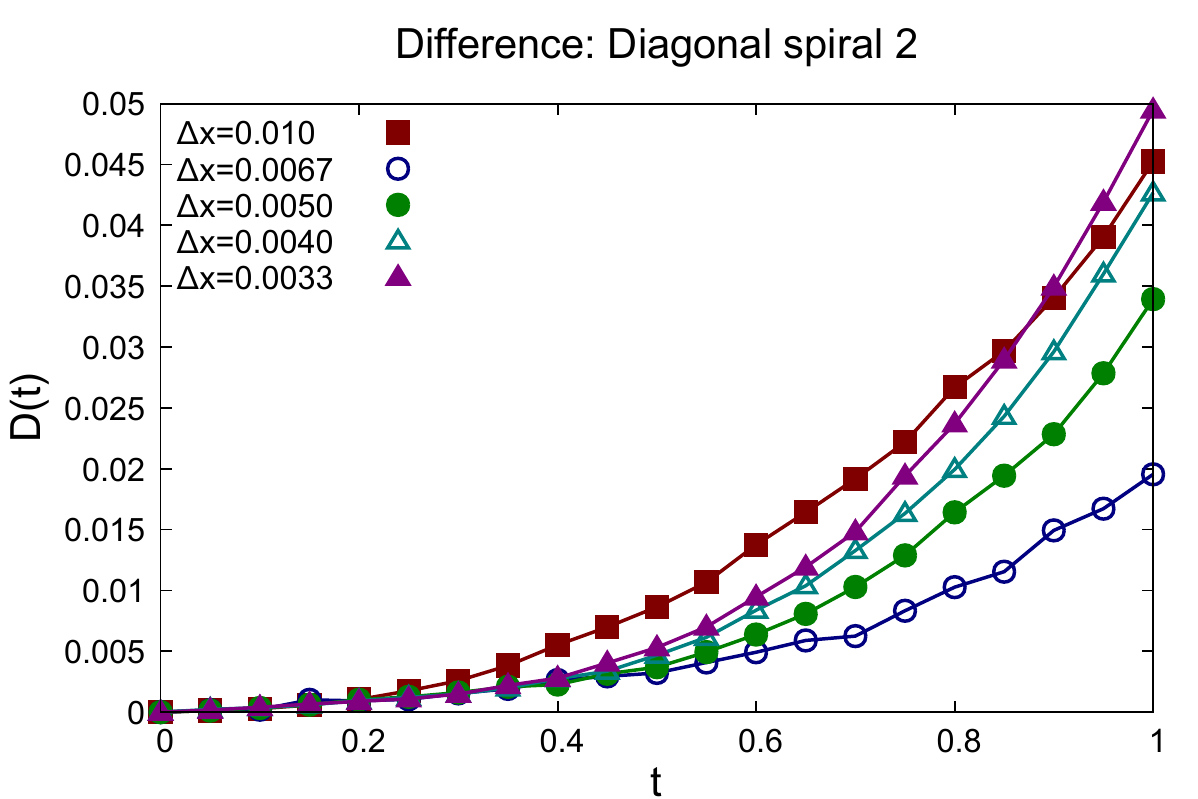}
  \caption{Graphs of functions $t \mapsto \mathcal{D}(t)$ for the
  diagonal spiral with a fixed center radius
  $\rho = 0.02-10^{-8}$(left), 
  and with a reduced center radius
  $\rho = (4 - 10^{-8}) \Delta x$(right).
  }
  \label{diff diag}
 \end{center}
\end{figure}

The left figure of
Figure \ref{diff diag}
is a graphs of $\mathcal{D} (t)$ for
$s=2,3,4,5,6$ with a fixed center radius
$\rho = 0.02 - 10^{-8}$.
One can find that $\mathcal{D} (t)$ is reduced
by choosing smaller $\Delta x$,
and the smallest $\mathcal{D} (t)$ is
the case $\Delta x = 0.0033$ ($s=6$).
Our numerical simulations show that
the differences
are less than 4\% of $|W|$ if $s \ge 4$.
Note that $\rho \approx 4 \Delta x$ when $s=4$.

Because of the above results, 
we choose $\rho \approx 4 \Delta x$ for accurate simulations
with a reduced center radius $\rho = O(\Delta x)$.
The right figure of
Figure \ref{diff diag} presents graphs of $\mathcal{D} (t)$
for $s = 2, 3, 4, 5, 6$ with $\rho = (4 - 10^{-8}) \Delta x$.
One can find that the differences are less than
$5 \%$ for all cases, although the worst one is
that with $\Delta x = 0.0033$ ($s=6$).
Note that the cases of $\Delta x=0.0050$ ($s=4$)
in both figures of Figure \ref{diff diag} are the same.

\subsection{Triangle spiral}
\label{sec: triangle}

\begin{figure}[htbp]
 \begin{center}
  \includegraphics[scale=0.55]{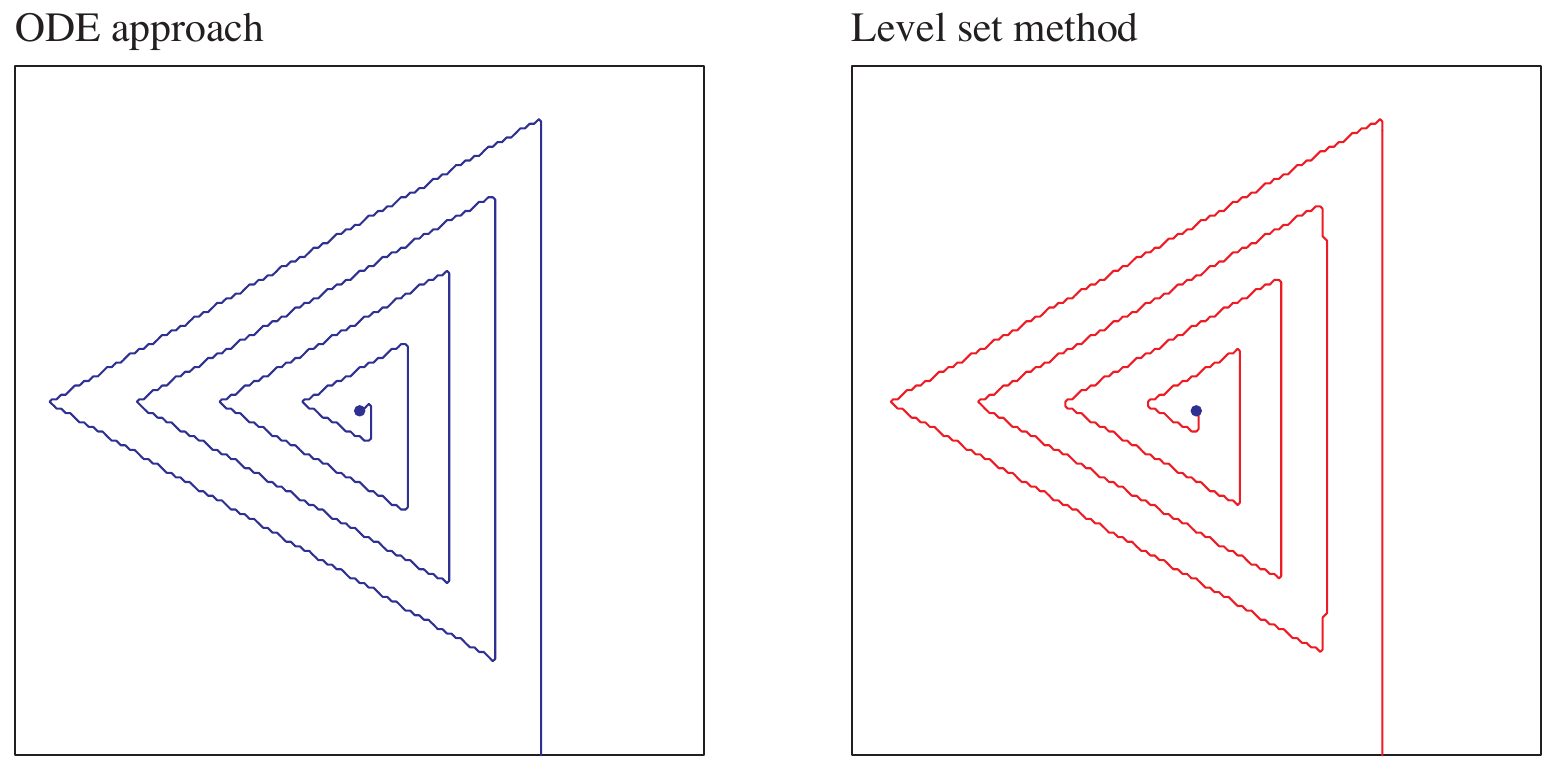}
  \caption{Profiles of the triangle spiral at $t=0.8$.
  The level set method is calculated with
  $\rho = 0.02 - 10^{-8}$ and $\Delta x = 0.0050$.}
  \label{profile tri2}
 \end{center}
\end{figure}

Finally, we examine a triangle spiral
as an asymmetric case of $\gamma$ or $\gamma^\circ$.
To give its settings, we first give $\gamma^\circ$.
Because of the normalizing assumption \eqref{normalization},
we set
\begin{align*}
 \gamma^\circ (p) = \max_{0 \le j \le 2} m_j \cdot \pi 
 \quad \mbox{with} \
 m_j = \left( \cos \frac{2 \pi j}{3}, \sin \frac{2 \pi j}{3} \right).
\end{align*}
Then, $\mathcal{W}_\gamma = \{ p \in \mathbb{R}^2 ; \ 
\gamma^\circ (p) \le 1 \}$ implies that
\[
 \varphi_j = \frac{2 \pi j}{3}, \quad \ell_j = 2 \sqrt{3}
\]
since $\mathcal{W}_\gamma$ is an equilateral triangle
whose vertices are at $(1, \pm \sqrt{3})$ and $(-2,0)$.
On the other hand, from the computation 
as in \S \ref{sec: how to define density} we obtain
\begin{align*}
 & \gamma (p) = \max_{0 \le j \le 2} n_j \cdot p \quad
 \mbox{with} \ n_j = 2 \left( \cos \frac{(2j+1) \pi}{3},
 \sin \frac{(2j+1) \pi}{3} \right), \\ 
 & \mbox{and then} \quad
 \tilde{\gamma} (p)
 = \gamma (- p)
 = \max_{0 \le j \le 2} \tilde{n}_j \cdot p \quad
 \mbox{with} \ \tilde{n}_j = 2 \left( \cos \frac{2 \pi j}{3},
 \sin \frac{2 \pi j}{3} \right).
\end{align*}
Note that $Q_j$ in \eqref{char form of gamma} is given as
\begin{align*}
 & Q_j = \{ p \in \mathbb{R}^2 ; \ g_j (p) \ge g_k (p) \ \mbox{for} \ 
 k \neq j \}  
 = \{ p \in \mathbb{R}^2 ; \ 
 \min_{k \neq j} (g_j (p) - g_k (p)) \ge 0 \} \\
 & \qquad \mbox{with} \quad
 g_0 (p) = 2 p_1, \ g_1 (p) = - p_1 + \sqrt{3} p_2, \
 g_2 (p) = - p_1 + \sqrt{3} p_2.
\end{align*}
Then, we obtain
\[
 \tilde{\gamma} (p)
 \approx \sum_{j=0}^2 (\tilde{n}_j \cdot p)
 \zeta (f_j (p); 1, 0),
 \quad
 \tilde{\xi} (p)
 \approx \sum_{j=0}^2  
 \zeta (f_j (p); 1, 0) \tilde{n}_j,
\]
where
$f_j (p) = \min_{k \neq j} (g_j (p) - g_k (p))$.
We calculate \eqref{ODE system1}--\eqref{ODE system2}
or \eqref{lv mcf}--\eqref{lv nbc} for the evolution equation
\[
 V_\gamma = 1 - 0.01 H_\gamma
 \quad (\beta \equiv 1, \ U=1, \ \rho_c = 0.01)
\]
on the time interval $[0,0.8]$ with the above anisotropic setting.
This time interval is chosen so that $\Gamma_D (t)$ does not
touch to $\partial \Omega$ for $t \in [0,0.8]$.
Figure \ref{profile tri2} are profiles of the diagonal spiral
at $t=0.8$ with the above setting.

\begin{figure}[htbp]
 \begin{center}
  \includegraphics[scale=0.5]{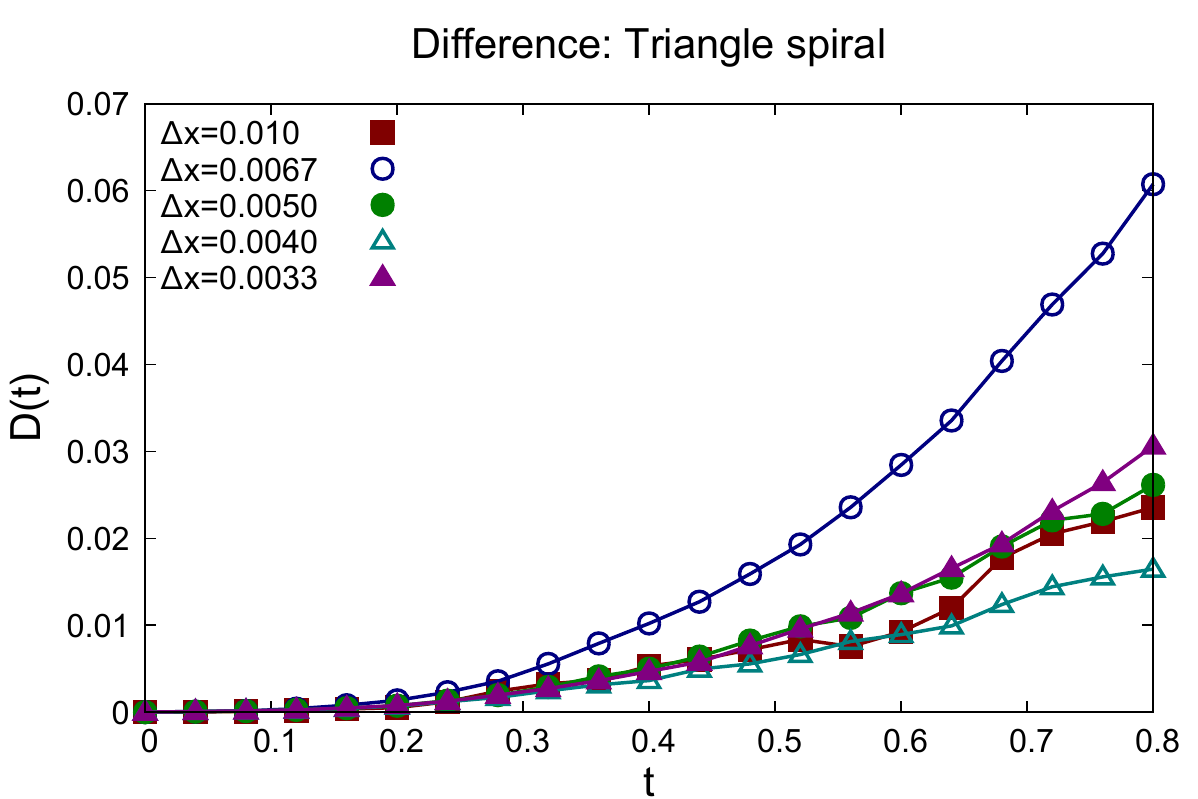}
  \includegraphics[scale=0.5]{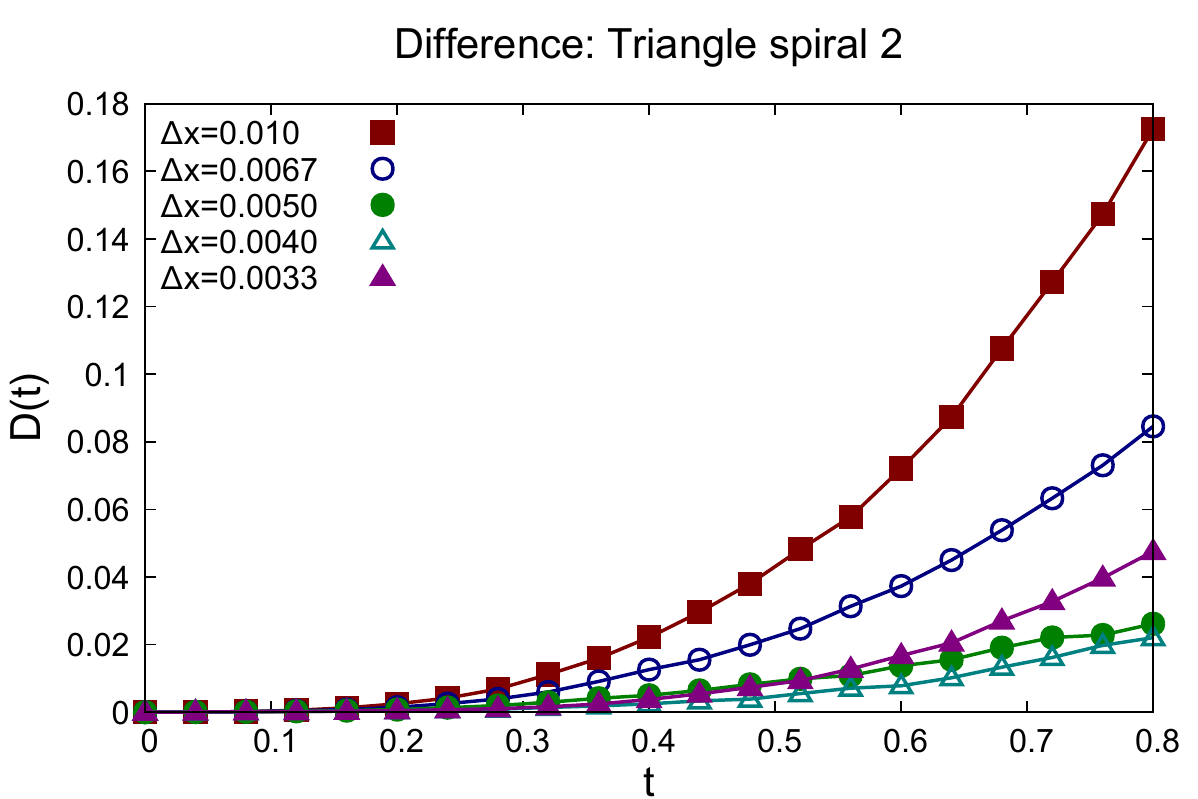}
  \caption{Graphs of functions $t \mapsto \mathcal{D} (t)$
  for the triangle spiral with a fixed center radius
  $\rho = 0.02 - 10^{-8}$(left),
  and with a reduced center radius
  $\rho = (4 - 10^{-8}) \Delta x$(right).
  }
  \label{diff tri}
 \end{center}
\end{figure}

The left figure of 
Figure \ref{diff tri} presents graphs of $\mathcal{D} (t)$
with $s=2, 3, 4, 5, 6$ with a fixed center radius 
$\rho = 0.02 - 10^{-8}$.
One can find that the differences are less than
$4 \%$ for the all cases except $\Delta x = 0.0067$ ($s=3$), 
and the best one is
that with $\Delta x = 0.0040$ ($s=5$).

From the analogy of the diagonal spiral case,
we choose $\rho \approx 4 \Delta x$ 
as a reducing center radius $\rho = O(\Delta x)$.
The right figure of
Figure \ref{diff tri} presents graphs of $\mathcal{D} (t)$
with $s=2, 3, 4, 5, 6$ with a fixed center radius 
$\rho = (4 - 10^{-8}) \Delta x$.
One can find that the differences are less than
$5 \%$ when $\Delta x \le 0.0050$ ($s \ge 4$).
Note that the cases of $\Delta x=0.0050$($s=4$) in both figure of
Figure \ref{diff tri} are the same.

\section{Conclusion}

In this paper, we compared the discrete model
as in \cite{IO:DCDS-B} and the level set method
as in \cite{OTG:2015JSC} for evolving spirals
by the crystalline eikonal-curvature flow \eqref{geo mcf}.
Note that the level set equation includes
the derivative of a piecewise linear energy density function.
For this problem, 
we introduced an approximation of
level set equation for the crystalline curvature flow,
which is established with the approximation of 
the characteristic function as in \cite{ETT2005JCP}.
To measure the difference between the two curves
obtained by the discrete model and the level set method,
we introduced an area difference function defined
by \eqref{area difference}.
It is consist of the $L^1$ difference 
of the height function as in \cite{OTG:2015JSC}
with the step-height $h_0 = 1$.
Note that the discrete and level set models are
slightly different on the boundary condition
at the center and the outer boundary of the domain.
However, we found that the area differences of
these models are less than $5 \%$ of the area
of the domain for square, diagonal and triangle
spirals as in \S \ref{sec: numerical results}
when the resolution of the numerical lattice
is enough high and the radius of the center for the
level set method is suitably small.

\medskip

\end{document}